\documentclass[smallextended,referee,envcountsect]{svjour3}
\usepackage[ruled, vlined]{algorithm2e}
\usepackage{amsmath}
\usepackage{amssymb}
\usepackage{array}
\usepackage{booktabs}
\usepackage{float}
\usepackage[margin=1in]{geometry}
\usepackage{graphicx}
\usepackage{hyperref}
\usepackage{tabu}
\usepackage{tabularx, ragged2e}

\let\proof\relax

\usepackage{amsthm}

\def\blackslug{\hbox{\kern1pt\vrule height6pt width4pt  depth1pt\kern1pt}}

\def\edp{\penalty 500\hbox{\quad\blackslug}\ifmmode\else\par
    \vskip4.5pt plus3pt minus2pt\fi}
\def\proof{\par\noindent{\bf Proof.\enspace}\rm}
\def\qed{\penalty 500\hbox{\quad\blackslug}\ifmmode\else\par
    \vskip4.5pt plus3pt minus2pt\fi}

\begin{document}
\title{Computing skeletons for rectilinearly-convex obstacles in the rectilinear plane\thanks{This work was supported by an Australian Research Council Discovery Grant}}
\author{Marcus Volz \and Marcus Brazil \and Charl Ras \and Doreen Thomas}

\institute{Corresponding author: Marcus Volz
\at School of Engineering, The University of Melbourne, Victoria 3010, Australia\\\email{marcus.volz@unimelb.edu.au}\\Tel. +61 402 859 710
\and
Marcus Brazil \and Doreen Thomas
\at School of Engineering, The University of Melbourne, Victoria 3010, Australia
\and
Charl Ras
\at School of Mathematics and Statistics, The University of Melbourne, Victoria 3010, Australia
}

\maketitle

\begin{abstract}
  We introduce the concept of an obstacle \emph{skeleton} which is a set of line segments inside a polygonal obstacle $\omega$ that can be used in place of $\omega$ when performing intersection tests for obstacle-avoiding network problems in the plane. A skeleton can have significantly fewer line segments compared to the number of line segments in the boundary of the original obstacle, and therefore performing intersection tests on a skeleton (rather than the original obstacle) can significantly reduce the CPU time required by algorithms for computing solutions to obstacle-avoidance problems. A \emph{minimum skeleton} is a skeleton with the smallest possible number of line segments. We provide an exact $O(n^2)$ algorithm for computing minimum skeletons for rectilinear obstacles in the rectilinear plane that are rectilinearly-convex. We show that the number of edges in a minimum skeleton is generally very small compared to the number of edges in the boundary of the original obstacle, by performing experiments on random rectilinearly-convex obstacles with up to 1000 vertices.
  \subclass{Geometry \and Convex and discrete geometry \and Computer science}
\end{abstract}

\section{Introduction}

Obstacle-avoiding shortest network problems arise in many applications in industry. In these problems it may be necessary to perform intersection tests to determine whether a part of the network intersects one or more obstacles, since a network that intersects an obstacle is not a feasible solution. For polygonal obstacles, the number of edges can significantly contribute to the CPU time required by algorithms for computing shortest obstacle-avoiding networks. It is therefore desirable to reduce the number of edges that need to be considered for a given obstacle, since this will reduce the number of obstacle intersection tests that need to be performed.

Throughout this paper an \emph{obstacle} (denoted by $\omega$) is a simple polygon, i.e. a closed and bounded polygonal region that does not have holes and whose boundary does not intersect itself. We will make the further assumption that $\omega$ is in \emph{general position}, by which we mean that no three of its vertices are colinear. In obstacle-avoiding network problems the network is not permitted to intersect the interior of $\omega$. In this paper we introduce the concept of a \emph{skeleton}, which is a representation of an obstacle that consists of a set of line segments inside the obstacle. A skeleton can have a significantly reduced number of line segments compared to the number of line segments in the boundary of the original obstacle, and therefore performing intersection tests on a skeleton (rather than the original obstacle) can significantly reduce the CPU time required by algorithms for computing shortest obstacle-avoiding networks.

\begin{figure}[h]
\centering
\includegraphics[width=13cm]{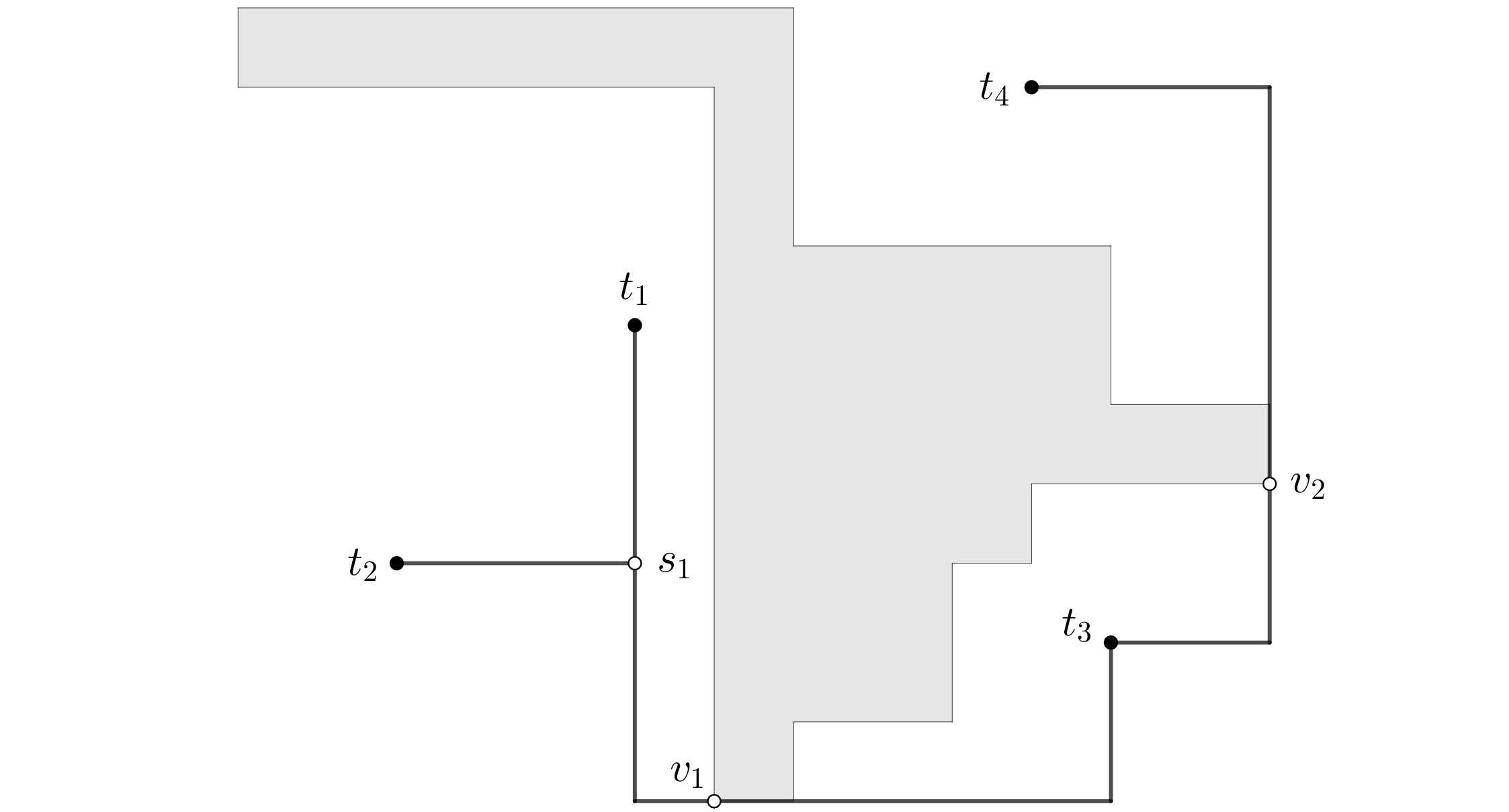}
\caption{A shortest obstacle-avoiding rectilinear network interconnecting points $t_1$, $t_2$, $t_3$ and $t_4$. The network has two straight edges ($t_1s_1$ and $t_2s_1$) and four bent edges ($s_1v_1, v_1t_3, t_3v_2$ and $v_2t_4$).}\label{fig-oarsmt-example}
\end{figure}

In this paper we focus on computing minimum skeletons for obstacles in the context of shortest networks with respect to the rectilinear metric (Figure~\ref{fig-oarsmt-example}). Obstacle-avoiding rectilinear network problems and shortest rectilinear path problems have been well studied (\cite{lee1996rectilinear}, \cite{macgregor1979steiner}, \cite{ganley1994routing}, \cite{lin2007efficient}, \cite{li2008obstacle}, \cite{lin2008obstacle}, \cite{long2008eboarst}, \cite{liu2009high}, \cite{huang2010obstacle}, \cite{huang2011exact}, \cite{ajwani2011foars}, \cite{huang2011construction}, \cite{liu2012obstacle}, \cite{huang2013obsteiner}, \cite{chow2014obstacle}, \cite{held2014fast}, \cite{brazil2015optimal}, \cite{volz2019}) and have a range of applications including VLSI design and motion planning. When represented as an embedding in the Euclidean plane, each edge is then a rectilinear shortest path between its endpoints; that is, a shortest path composed only of horizontal and vertical line segments.

Shortest rectilinear networks can be constructed from two types of edges (which represent geodesics in the rectilinear metric): a \emph{straight edge} is a single line segment that is either horizontal or vertical; and a \emph{bent edge} consists of a series of horizontal and vertical line segments where each adjacent pair of orthogonal line segments meet at a \emph{corner point}. Any bent edge connecting a given pair of points can be embedded with a single corner point and has exactly two such embeddings. In many problem contexts, it is sufficient to consider shortest networks for which each edge has at most one corner point (see for example the discussion below on obstacle-avoiding rectilinear Steiner trees).

On this basis we provide a formal definition of skeletons for obstacles in the context of rectilinear shortest networks as follows (note that a set of embedded edges $S$ is considered to be \emph{inside} $\omega$ if $\bigcup\limits_{s_i\in S} s_i\subseteq\omega$):

\begin{definition} Let $S$ be a set of closed line segments inside a polygonal obstacle $\omega$. Then $S$ is said to be a \emph{skeleton} for $\omega$ if for any given pair of points outside the interior of $\omega$, if every rectilinear shortest path with at most one corner point between the pair of points meets the interior of $\omega$, then each such path intersects an element of $S$. A \emph{minimum skeleton} $S^*$ is a skeleton with the smallest possible number of line segments.\label{def-skeleton}
\end{definition}

\begin{figure}[h]
\centering
\includegraphics[width=13cm]{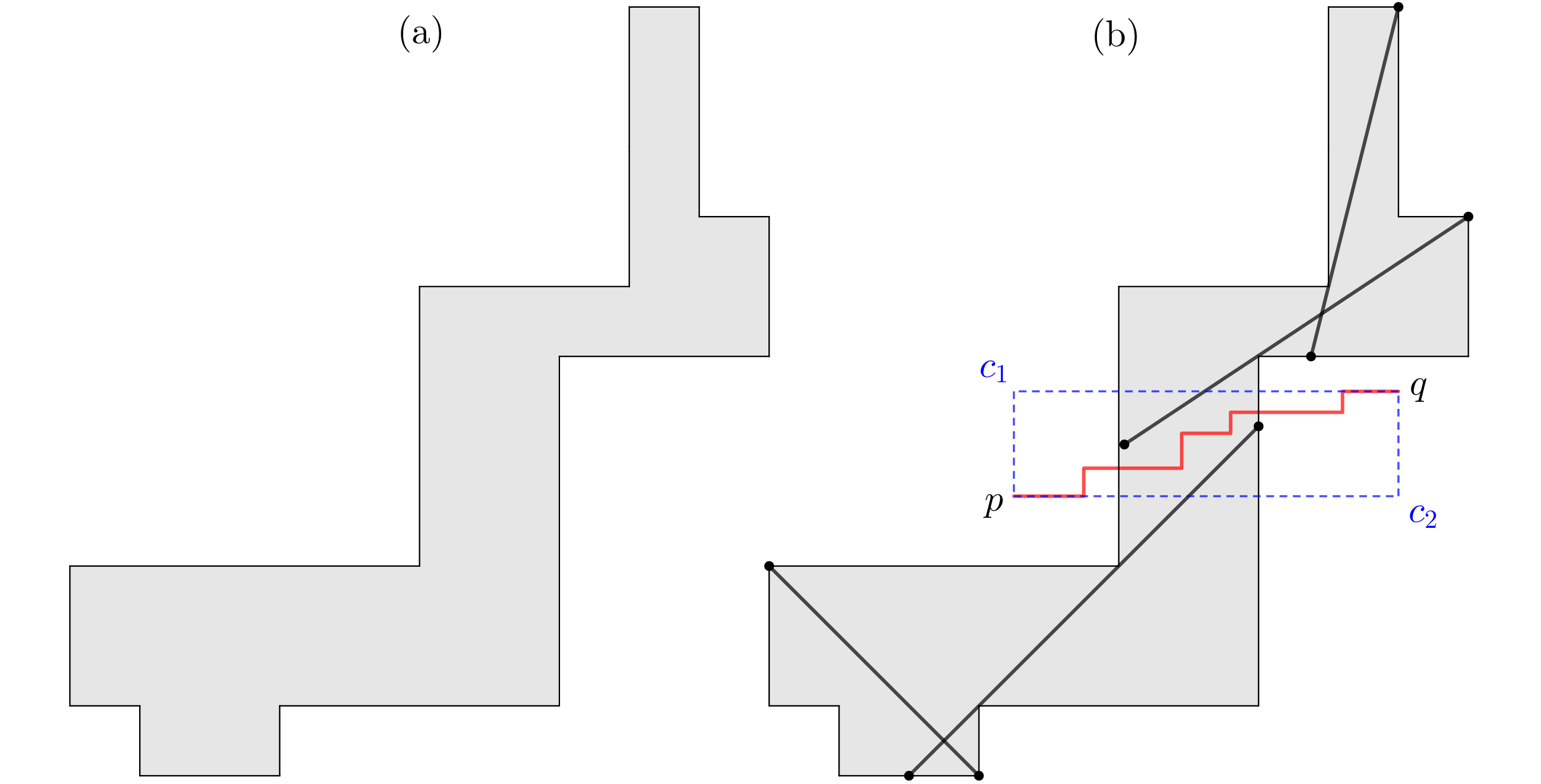}
\caption{(a) A rectilinearly-convex obstacle $\omega$. (b) A minimum skeleton for $\omega$ with four line segments.}\label{fig-skeleton-example}
\end{figure}

An example of a minimum skeleton for an obstacle $\omega$ is given in Figure~\ref{fig-skeleton-example}. Note that shortest rectilinear paths between points outside the interior of $\omega$ with more than one corner point (i.e. zigzag paths) that pass through $\omega$ need not intersect a skeleton. For example in Figure~\ref{fig-skeleton-example}, a shortest rectilinear path between $p$ and $q$ with multiple corner points (shown as a red line) intersects $\omega$ but does not intersect the skeleton. The skeleton is legitimate, however, because the two shortest rectilinear paths with one corner point (i.e. $pc_1q$ and $pc_2q$) intersect the skeleton.

The requirement in Definition~\ref{def-skeleton} that a skeleton edge be inside $\omega$ ensures that all shortest rectilinear paths with at most one corner point that intersect the interior of a skeleton also intersect $\omega$. This requirement makes it necessary to treat obstacles individually when determining their skeletons (since any one skeleton constructed for multiple disjoint obstacles would in part lie outside the interior of the obstacles). The related question, of whether a rectilinear path that intersects a skeleton edge (say, at an endpoint of the skeleton edge) actually enters the interior of the obstacle or simply runs along part of the obstacle boundary, is addressed in Section~\ref{sect-skeletons-oarsmt}.

In this paper we focus on computing minimum skeletons for obstacles that are themselves rectilinear, i.e. simple polygons for which each edge of the polygon is either horizontal or vertical. In particular, we focus on rectilinear obstacles that are \emph{rectilinearly-convex}, meaning that any two points in $\omega$ can be joined by a shortest rectilinear path that is inside $\omega$. This is a reasonable restriction of the problem to study initially, as it provides insight into the structure of minimum skeletons, and the intuition gained through studying minimum skeletons for rectilinearly-convex obstacles can then be extended to other contexts. Moreover, rectilinearly-convex obstacles are often found in VLSI problem instances (see for instance the many examples in the SteinLib database~\cite{koch2001steinlib}). For the remainder of this paper, the term \emph{rectilinearly-convex} will be used to refer to an obstacle that is both rectilinear and rectilinearly-convex.

\subsection{Obstacle-avoiding rectilinear Steiner trees and other applications of skeletons}\label{sect-skeletons-oarsmt}

We now discuss the application of skeletons to the \emph{obstacle-avoiding rectilinear Steiner tree problem}, a problem which is applicable to physical networks in VLSI design (or microchip design); see Section~3.6 of~\cite{brazil2015optimal}. The following properties of rectilinear Steiner trees are from Section 3.1 of~\cite{brazil2015optimal}. Given a set of simple rectilinearly-convex obstacles in the Euclidean plane with a collective set $V$ of vertices, called \emph{virtual terminals}, and a set $N$ of points, called \emph{terminals}, that are outside the interiors of the obstacles, an \emph{obstacle-avoiding rectilinear minimum Steiner tree (OARMST)} is a shortest network interconnecting the terminals, where each edge of the network is composed of horizontal and vertical line segments, and no point (vertex or interior point of an edge) on the network lies in the interior of an obstacle. This is equivalent to finding a shortest embedded obstacle-avoiding network in which the length of each edge is determined by the \emph{rectilinear metric} (i.e. the length of an edge is the length of a shortest rectilinear path between the endpoints of the edge).

An exact algorithm~\cite{huang2013obsteiner} exists for computing OARMSTs which is based on the \emph{GeoSteiner} algorithm~\cite{warme2000exact}. GeoSteiner uses the fact that OARMSTs decompose into \emph{full components}, each of which is a \emph{full} minimum Steiner tree (FST) on a subset of $N \cup V$, meaning it is a minimum Steiner tree in which every terminal and every virtual terminal of each component has degree one. In the \emph{generation phase} of GeoSteiner, a set of candidate full components is computed such that the set is guaranteed to contain all full components of a minimum solution. Then in the \emph{concatenation phase}, a minimum interconnection network is constructed from the candidate full components.

The generation phase of GeoSteiner consists of a bottom-up construction of full components via the construction of branches in which one of the edges is a ray. At each stage in the construction of the branch it is necessary to check for intersections with the boundary of any obstacle. If such an intersection occurs, the branch, and hence the family of full components it would have generated, can be discarded. Large complex obstacles with many edges can therefore significantly contribute to the CPU time required by the algorithm. It is therefore desirable to reduce the number of edges that need to be considered for a given obstacle, since this will reduce the number of obstacle intersection tests that need to be performed in the generation phase of GeoSteiner.

The most effective GeoSteiner algorithms assume that each edge of the minimum network contains at most one corner point. This also makes sense in VLSI applications as any change in direction in an edge incurs a cost, since it usually involves moving from one layer of the microchip to another. Other related problems with applications to VLSI physical design include obstacle avoiding minimum rectilinear spanning trees, and obstacle avoiding minimum rectilinear visibility graphs, where again in each case all edges are assumed to contain at most one corner point. In all of these cases, any rectilinear edge embedded with at most one corner point that does not intersect any skeleton edges has an embedding that is obstacle avoiding. It now remains to specify the conditions under which such a rectilinear edge that does intersect a skeleton edge also intersects the interior of the obstacle.

Let $s$ be a closed line segment inside a rectilinearly-convex obstacle $\omega$ with both endpoints on the boundary of $\omega$. We can classify each point on $s$ as either a \emph{boundary point} if it lies on the boundary of $\omega$, or a non-boundary point if it lies in the interior of $\omega$. Let $P$ be a rectilinear path connecting two points outside the interior of $\omega$, and containing at most one corner point. We first define what it means for $s$ to \emph{block} $P$.

\begin{definition}  With $s$, $P$ and $\omega$ defined as above, we say that $s$ \emph{blocks} $P$ if:
\begin{enumerate}
  \item $P$ intersects $s$ at a non-boundary point; or
  \item $P$ meets a boundary point of $s$ not at a corner point of $P$, \emph{and} the direction of $P$ at this meeting point is not equal to the direction of any edge of the boundary of $\omega$ containing this boundary point of $s$; or
  \item $s$ is neither horizontal nor vertical, \emph{and} $P$ meets an endpoint of $s$ at a corner point of $P$, \emph{and} $P$ includes part or all of the  interior of a rectilinear bounding box edge of $s$ that intersects the interior of $\omega$.
\end{enumerate}
\end{definition}

It is straightforward to see that this definition is precisely what is required to get the following theorem.

\begin{theorem} Let $S$ be a skeleton for a rectilinearly-convex obstacle $\omega$. Let $P$ be a rectilinear path connecting two points outside the interior of $\omega$, and containing at most one corner point. Then an edge of $S$ blocks $P$ if and only if every embedding of $P$ with at most one corner point intersects the interior of $\omega$. \end{theorem}

Once a skeleton has been constructed for a given obstacle $\omega$, information about how each skeleton edge interacts with the boundary of $\omega$ can be recorded. This means that it is simple to check whether a skeleton edge blocks a given rectilinear path using only this supplementary information, and without needing to refer back to $\omega$.

\section{General properties of rectilinearly-convex obstacles}

A \emph{rectilinearly-convex obstacle} $\omega$ is a simple polygon for which each edge of the polygon is either horizontal or vertical and such that for any two points in $\omega$ there is a shortest rectilinear path between the two points that is inside $\omega$. In this section we identify some useful general properties of rectilinearly-convex obstacles.

\subsection{Classification of rectilinearly-convex obstacles}\label{sect-classification}

Throughout this paper we will use the terms ``edges of $\omega$" and ``vertices of $\omega$" as shorthand to refer to the edges and vertices of the boundary of $\omega$ (since $\omega$ is a region). The \emph{extreme edges} of $\omega$ are the edges of $\omega$ that intersect its \emph{bounding box}, where the bounding box is the smallest closed axis-oriented rectangle that encloses $\omega$. The term \emph{bounding box} will be similarly applied more generally to any set of line segments, such as skeletons. The function that returns the bounding box for a given set of line segments will be denoted by $B(\cdot)$. An \emph{extreme corner} of $\omega$ is a vertex of $\omega$ at the intersection of two extreme edges. A rectilinearly-convex obstacle has exactly four extreme edges and up to four extreme corners.

A \emph{staircase walk} of $\omega$ is a shortest rectilinear path on the boundary of $\omega$ between adjacent extreme edges (including the extreme edges themselves). A pair of extreme corners (edges) will be called \emph{adjacent} if they coincide with adjacent vertices (edges) of the bounding box of $\omega$, or \emph{opposite} otherwise. An extreme corner $c$ and an extreme edge $e$ are \emph{opposite} if $e$ lies on an edge of the bounding box of $\omega$ that does not intersect $c$. For example, if $c$ lies at the bottom-left corner of $B$, then the right and top extreme edges of $\omega$ are opposite to $c$.

Rectilinearly-convex obstacles can be categorised into six types based on the number and relative locations of extreme corners (see Figure~\ref{figS}):

\begin{figure}[h]
\centering
\includegraphics[width=\textwidth]{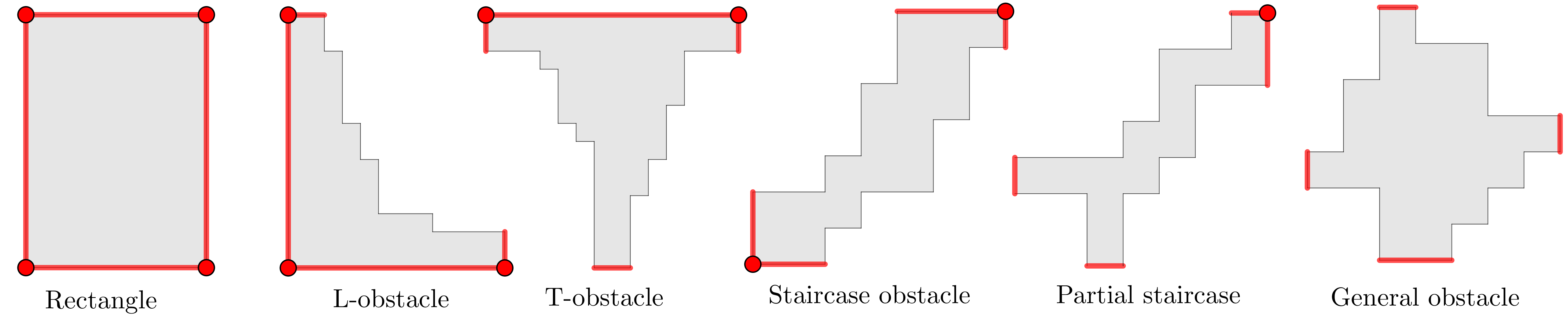}
\caption{Classification of rectilinearly-convex obstacles based on the relative locations of extreme corners. Extreme edges are shown as red lines, and extreme corners as red dots.}\label{figS}
\end{figure}

\begin{itemize}
  \item A \emph{rectangle} has exactly four extreme corners.
  \item An \emph{L-obstacle} has exactly three extreme corners.
  \item A \emph{T-obstacle} has exactly two extreme corners, which are adjacent.
  \item A \emph{staircase obstacle} has exactly two extreme corners, which are opposite.
  \item A \emph{partial staircase} has exactly one extreme corner.
  \item A \emph{general obstacle} has no extreme corners.
\end{itemize}

\subsection{Sub-classification of general obstacles}\label{sect-general-subtypes}

Let $\omega$ be a general obstacle. A pair of parallel extreme edges of $\omega$ will be said to \emph{overlap} if the orthogonal projection of one edge onto the other is not empty, and a pair of non-overlapping extreme edges will be called \emph{positively (negatively) sloped} if the absolute gradient of the line segment between the midpoints of the two edges is positive (negative). Up to symmetry, $\omega$ can be sub-classified into four types as shown in Figure~\ref{fig-general-obstacles}:

\begin{figure}[h]
\centering
\includegraphics[width=14cm]{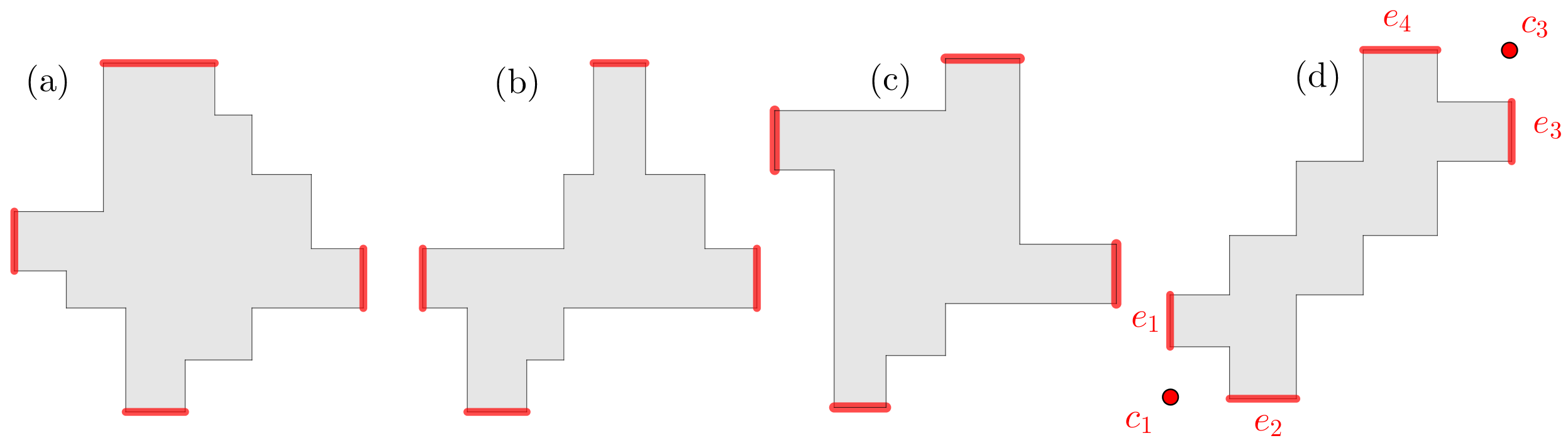}
\caption{Sub-classification of general obstacles based on the relative positions of parallel extreme edges.}\label{fig-general-obstacles}
\end{figure}

\begin{itemize}
  \item Type (a): Both pairs of parallel extreme edges overlap.
  \item Type (b): Exactly one pair of parallel extreme edges overlaps.
  \item Type (c): Neither pair of parallel extreme edges overlap and the two pairs of parallel extreme edges have different slopes.
  \item Type (d): Neither pair of parallel extreme edges overlap and the two pairs of parallel extreme edges have the same slope.
\end{itemize}

\subsection{Obstacle ends}

Let $e_i$ and $e_j$ be adjacent extreme edges of a rectilinearly-convex obstacle $\omega$ and assume that $e_i$ is horizontal and $e_j$ is vertical. Then $\{e_i, e_j\}$ is called an \emph{end} of $\omega$ if: (1) $e_i$ and $e_j$ do not overlap with their corresponding opposite parallel extreme edges, and (2) there exists a corner $c$ of $B(\omega)$ such that $e_i$ is the closest horizontal extreme edge to $c$ in the horizontal direction and $e_j$ is the closest vertical extreme edge to $c$ in the vertical direction.

A staircase or partial staircase obstacle with overlapping parallel extreme edges has no ends; otherwise, a staircase obstacle has two ends, one corresponding to each extreme corner, while a partial staircase has one end corresponding to its extreme corner $c$ and the other end corresponding to the two extreme edges that are not incident to $c$. Rectangles, L-obstacles and T-obstacles do not have ends, since their pairs of parallel extreme edges overlap. Type (a), (b) and (c) general obstacles do not have ends, while Type (d) general obstacles have two ends. For example, the obstacle in Figure~\ref{fig-general-obstacles} (d) has two ends: $\{e_1,e_2\}$ and $\{e_3,e_4\}$ corresponding to $c_1$ and $c_3$.

\subsection{Point and edge visibility}

Two points, $p$ and $q$, inside a rectilinearly-convex obstacle $\omega$ will be called \emph{mutually visible} if the line segment between $p$ and $q$ is inside $\omega$. Two line segments $s_1$ and $s_2$ in $\omega$ are mutually visible if there exists a point $p_1$ on $s_1$ and a point $p_2$ on $s_2$ such that $p_1$ and $p_2$ are mutually visible. A point $p$ and a line segment $s$ in $\omega$ are mutually visible if there exists a point $q$ on $s$ such that $p$ and $q$ are mutually visible.

A \emph{visibility edge} $e$ is any line segment that is inside $\omega$. An \emph{auxiliary point} is an endpoint of a visibility edge that lies in the interior of an edge in the boundary of $\omega$. We define the following special types of visibility edge configurations, which may or may not exist for a given obstacle (refer to Figure~\ref{fig-vis-edges}):

\begin{figure}[h]
\centering
\includegraphics[width=\textwidth]{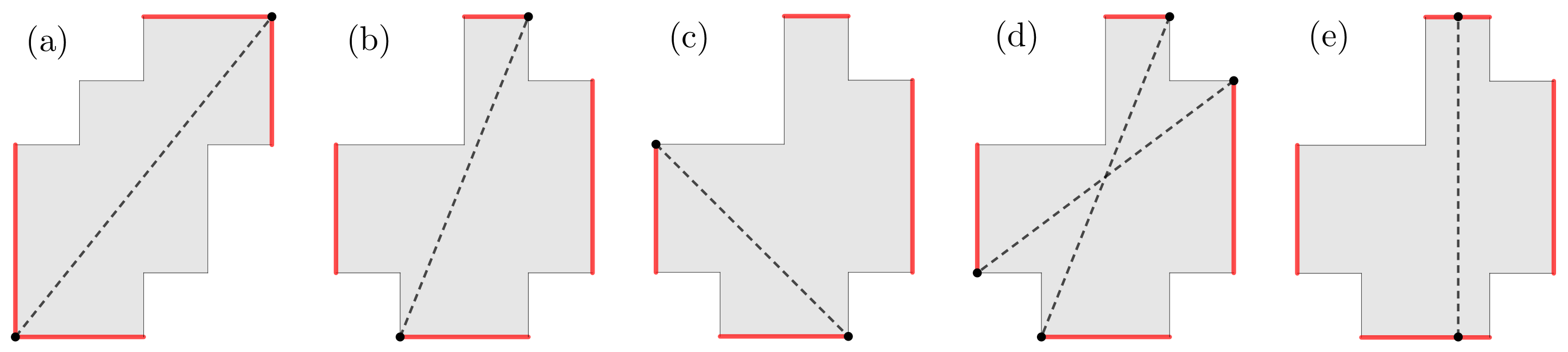}
\caption{Special types of edge configurations: (a) Diagonal. (b) Opposite extreme visibility edge. (c) Maximum length adjacent extreme visibility edge. (d) Cross. (e) Perpendicular extreme visibility edge.}\label{fig-vis-edges}
\end{figure}

\begin{itemize}
  \item A \emph{diagonal} of $\omega$ is a line segment inside $\omega$ between a pair of opposite extreme corners.
  \item An \emph{opposite (adjacent) extreme visibility edge} is a line segment inside $\omega$ that has endpoints on parallel (orthogonal) extreme edges.
  \item A \emph{maximum length adjacent extreme visibility edge} is an adjacent extreme visibility edge with the maximum length among all adjacent extreme visibility edges for a given pair of extreme edges (for example, if $e_1$ and $e_2$ are the left and bottom extreme edges respectively, then the maximum length extreme visibility edge connects the top endpoint of $e_1$ and the right endpoint of $e_2$).
  \item A \emph{perpendicular extreme visibility edge} is a visibility edge that is perpendicular to its corresponding extreme edge.
  \item A \emph{cross} is a pair $s_H,s_V$ of opposite extreme visibility edges, where $s_H$ ($s_V$) has a endpoint on each horizontal (vertical) extreme edge.
\end{itemize}

\section{General properties of skeletons}

In this section we present some general properties of skeletons. Consider a visibility edge $s$ inside a rectilinearly-convex obstacle. Then $s$ will be called a \emph{maximum length visibility edge} if it extends as far as possible in both directions to the boundary of $\omega$, subject to $s$ remaining inside $\omega$.

\begin{lemma} Let $\omega$ be a rectilinearly-convex obstacle. Then there exists a minimum skeleton $S^*$ for $\omega$ such that each edge in $S^*$ is a maximum length visibility edge.\label{lem-extend-edges}\end{lemma}

\proof Suppose $s$ is an edge of a skeleton $S^*$. Then, since $s$ is inside some obstacle $\omega$, any transformation of $S^*$ formed by extending $s$ at one or both ends to the boundary of $\omega$ is also a skeleton. The lemma immediately follows.\qed

For the remainder of this paper we will assume that all skeleton edges and visibility edges have the property of Lemma~\ref{lem-extend-edges}.

\begin{lemma} Let $\omega$ be a rectilinearly-convex obstacle and let $S^*$ be a skeleton for $\omega$. Then for each extreme edge $e_i, 1\leq i\leq 4$ of $\omega$, there exists an edge in $S^*$ with an endpoint on $e_i$.\label{lem-extreme-edge}\end{lemma}

\proof Without loss of generality, assume that there exists a skeleton that does not have an edge with an endpoint on the left extreme edge $e_L$ of $\omega$, and let $x(e_L)$ denote the $x$-coordinate of $e_L$. Since edges of the skeleton are closed line segments, there exists an $\epsilon > 0$ (where $\epsilon$ is strictly less than the length of the shortest horizontal edge of $\omega$) such that the $x$-distance of every point on the skeleton from $e_L$ is strictly greater than $\epsilon$. Let $l$ be a vertical line whose $x$-coordinate is $x(e_L) + \epsilon$. Then $l$ passes through $\omega$ without intersecting $S^*$, giving a contradiction.\qed

We define an \emph{extreme skeleton edge} to be a skeleton edge with at least one endpoint on an extreme edge.

\subsection{Horizontal and vertical projections}

The following is a useful property of skeletons.

\begin{lemma} Let $\omega$ be a rectilinearly-convex obstacle. If $S$ is a skeleton for $\omega$, then the horizontal and vertical projections of $S$ cover the horizontal and vertical projections of $\omega$.\label{lem-projections}\end{lemma}

\proof Without loss of generality, suppose that a subset $\Delta x$ of the horizontal projection of $\omega$ is not covered by the horizontal projection of $S$, and let $l$ be a vertical line that intersects the interior of $\Delta x$. Then $l$ passes through $\omega$ without intersecting $S$, and $S$ is therefore not a skeleton.\qed

\begin{figure}[h]
\centering
\includegraphics[width=10cm]{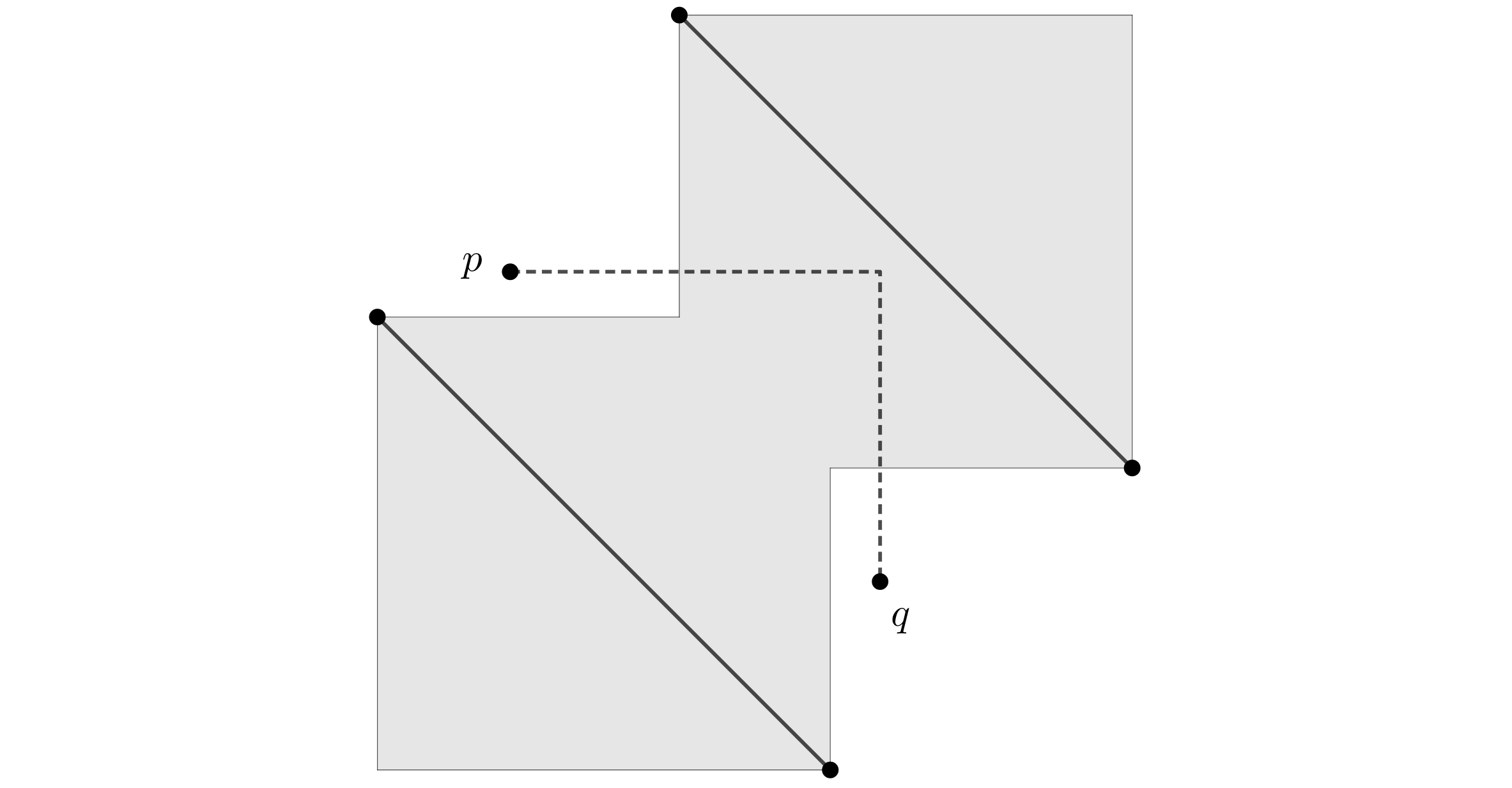}
\caption{The set of line segments is not a skeleton even though its horizontal and vertical projections cover the horizontal and vertical projections of the obstacle.}\label{fig-projections-covered-but-not-a-skeleton}
\end{figure}

Lemma~\ref{lem-projections} gives a necessary but not sufficient condition for $S$ to be a skeleton, as demonstrated by the example in Figure \ref{fig-projections-covered-but-not-a-skeleton}. In this example, the projections of the two skeleton edges cover the projections of the obstacle, however the two line segments clearly do not constitute a skeleton.

\subsection{Connectivity}

Let $S$ be a set of line segments embedded in the plane. Then we say $S$ is \emph{connected} if $S$ contains a path between every pair of points in $S$. The following lemma provides sufficient conditions for a set of line segments to be a skeleton.

\begin{lemma}
  Let $\omega$ be a rectilinearly-convex obstacle, and let $S$ be a set of line segments inside $\omega$ such that (1) $S$ intersects the four extreme edges of $\omega$ and (2) $S$ is connected. Then $S$ is a skeleton for $\omega$.\label{lem-connected}
\end{lemma}

\proof Let $S$ be a connected set of line segments inside $\omega$ that intersects the four extreme edges of $\omega$, and suppose that $S$ is not a skeleton for $\omega$. Then there exist points $p$ and $q$ outside the interior of $\omega$ such that all shortest rectilinear paths between $p$ and $q$ with at most one corner point intersect the interior of $\omega$, but at least one such path does not intersect $S$.

If $p$ and $q$ lie on a horizontal or vertical line, then there is a unique shortest rectilinear path between $p$ and $q$ (i.e. the line segment $pq$) which enters and exits $\omega$ at distinct points on the boundary of $\omega$. The obstacle $\omega$ can be partitioned into two regions, one on each side of $pq$, and each region contains exactly one extreme edge of $\omega$ that is parallel to $pq$, each of which is intersected by $S$. If $pq$ does not intersect $S$, it follows that $S$ must be disconnected, giving a contradiction.

Now suppose that $p$ and $q$ do not lie on a horizontal or vertical line. Then there are two shortest rectilinear paths with at most one corner point between $p$ and $q$, which we denote by $pc_1q$ and $pc_2q$. Since both $pc_1q$ and $pc_2q$ intersect $\omega$, it follows that $pc_1q$ cannot enter and exit $\omega$ on the same staircase walk; otherwise $pc_2q$ would be outside the interior of $\omega$ (due to the convexity of $\omega$).

Suppose that $pc_1q$ does not intersect $S$. If $c_1$ lies outside the interior of $\omega$, then at least one of the line segments $pc_1$ and $c_1q$ enters and exits $\omega$ at distinct points on the boundary of $\omega$, and the argument above (for the case where $p$ and $q$ are on a vertical or horizontal line) can be applied to arrive at the same contradiction. Otherwise, $c_1$ is in the interior of $\omega$, and $pc_1q$ enters and exits $\omega$ at two locations on distinct staircase walks. Therefore $pc_1q$ partitions $\omega$ into two regions, one on each side of $pc_1q$, and each region contains at least one extreme edge of $\omega$ (due to the convexity of $\omega$), each of which is intersected by $S$. If $pc_1q$ does not intersect $S$, it follows that $S$ must be disconnected, giving a contradiction.\qed

\subsection{Weak connectivity}

The converse of Lemma~\ref{lem-connected} is not true. That is, a skeleton is not necessarily connected, as demonstrated by the example in Figure \ref{fig-staircase-obstacle}, in which the three solid line segments form a minimum skeleton for the staircase obstacle. We therefore require a different kind of connectivity, which we will refer to as \emph{weak connectivity}, which is formally defined below:

\begin{figure}[h]
\centering
\includegraphics[width=10cm]{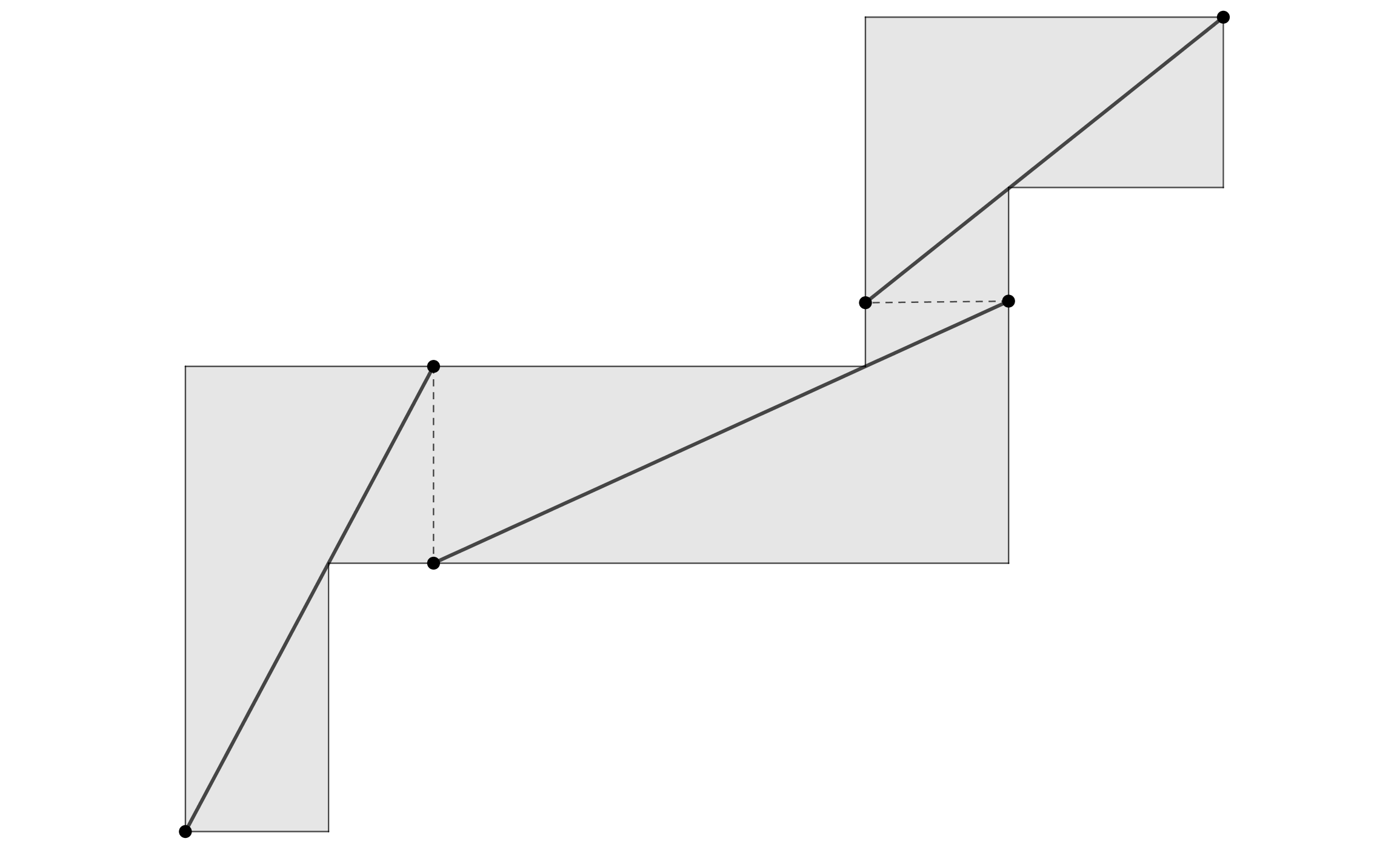}
\caption{An example of a minimum skeleton that is not connected.}\label{fig-staircase-obstacle}
\end{figure}

\begin{definition}
  Let $\omega$ be a rectilinearly-convex obstacle, and let $S$ be a set of line segments inside $\omega$. Let $\{S_i\}$ be the set of maximal connected sub-components of $S$.
  \begin{itemize}
    \item Two sub-components $S_j$ and $S_k$ will be called \emph{weakly connected} if $S_j\cap B_k\neq\emptyset$ and $S_k\cap B_j\neq\emptyset$ where $B_j$ and $B_k$ are the respective (closed) bounding boxes of $S_j$ and $S_k$.
    \item A skeleton is called \emph{weakly connected} if, when each sub-component is treated as a single vertex and an edge inserted between each pair of weakly connected sub-components, the resulting graph is connected.
  \end{itemize}
\end{definition}

We can now state necessary and sufficient conditions for a set of line segments to be a skeleton.

\begin{theorem}
   Let $\omega$ be a rectilinearly-convex obstacle and let $S$ be a set of line segments inside $\omega$. Then $S$ is a skeleton for $\omega$ iff (1) $S$ intersects the four extreme edges of $\omega$ and (2) $S$ is weakly connected.\label{thm-weakly-connected}
\end{theorem}

\begin{figure}[h]
\centering
\includegraphics[width=13cm]{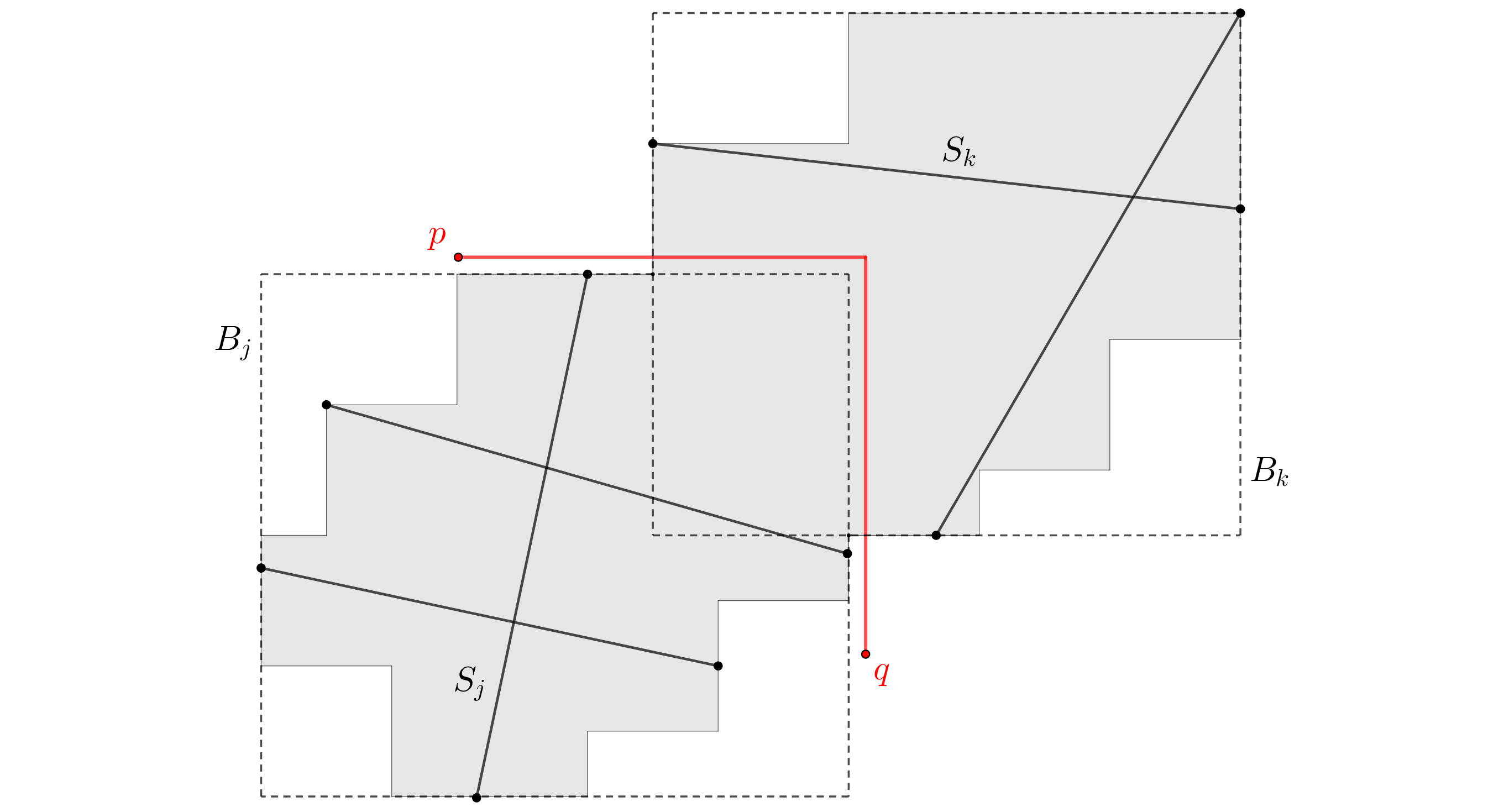}
\caption{Proof of Theorem~\ref{thm-weakly-connected}.}\label{fig-weakly-connected}
\end{figure}

\proof $(\rightarrow)$ Refer to Figure~\ref{fig-weakly-connected} and assume that $S$ is a skeleton. By Lemma~\ref{lem-extreme-edge}, $S$ necessarily intersects the four extreme edges of $\omega$. Now assume contrary to the theorem that $S$ is not weakly connected. Then $S$ can be partitioned into two disjoint subsets $S_j$ and $S_k$ such that $S_j\cap B_k = \emptyset$. Since $S$ is a skeleton, then by Lemma~\ref{lem-projections} the projection of $S_j\cup S_k$ must cover the projection of $\omega$, and therefore $B_j\cap B_k\neq\emptyset$. Hence there exists a shortest rectilinear path with one corner point between $p$ and $q$ (where the path follows the boundary of $B_j$ at some small distance outside $B_j$) that is not intersected by $S$ (as illustrated in Figure~\ref{fig-weakly-connected}). Therefore $S$ is not a skeleton, giving a contradiction.

$(\leftarrow)$ Now assume that $S$ intersects the four extreme edges of $\omega$ and that $S$ is weakly connected, and assume contrary to the theorem that $S$ is not a skeleton for $\omega$. Then there exist points $p$ and $q$ outside the interior of $\omega$ such that each shortest rectilinear path between $p$ and $q$ with at most one corner point intersects $\omega$ but at least one of these paths, say $pcq$, does not intersect $S$. Then $pcq$ partitions $S$ into two sets (as in the proof of Lemma~\ref{lem-connected}), and since each set necessarily intersects at least one extreme edge, part of the skeleton must lie strictly on each side of $pcq$. This contradicts the assumption that $S$ is weakly connected.\qed

\subsection{Uniqueness}

Figure~\ref{fig-staircase-obstacle} illustrates a case where the minimum skeleton is unique. We will see, however, that minimum skeletons are, in general, not unique. The example also demonstrates that both endpoints of an edge in a minimum skeleton can be auxiliary points.

\section{Minimum skeletons by obstacle type}

In this section we present methods for constructing minimum skeletons by obstacle type (based on the classification in Section~\ref{sect-classification}).

\subsection{Minimum skeletons for rectangles, L-obstacles and T-obstacles}

Let $|S|$ denote the cardinality of $S$ as a set of line segments.

\begin{lemma}
  Let $\omega$ be a rectilinearly-convex obstacle and let $S^*$ be a minimum skeleton for $\omega$. Then $|S^*| = 1$ iff $\omega$ has a diagonal.\label{lem-diagonal}
\end{lemma}

\proof $(\rightarrow)$ If $S^* = 1$, then the single edge $s$ in $S^*$ must intersect the four extreme edges of $\omega$ by Lemma~\ref{lem-extreme-edge}. This is only possible if the endpoints of $s$ are mutually visible opposite extreme corners of $\omega$. $(\leftarrow)$ A diagonal $\{d\}$ is a connected set of edges inside $\omega$ that intersects the four extreme edges of $\omega$. Therefore $\{d\}$ is a skeleton for $\omega$ by Lemma~\ref{lem-connected}. It is also a minimum skeleton since a skeleton clearly must have at least one edge.\qed

\begin{corollary}
  Let $\omega$ be a rectangle. Then $\{d\}$ is a minimum skeleton for $\omega$, where $\{d\}$ is either of the two diagonals of $\omega$.\label{cor-rectangle}
\end{corollary}

Corollary~\ref{cor-rectangle} is known in the literature, as discussed in Section 4.2.3 of~\cite{brazil2015optimal} (see Theorem 4.18). The result is also related to the observation in~\cite{huang2010obstacle} that for a given rectangular obstacle, only the two endpoints of a diagonal of the obstacle need to be considered when constructing an obstacle-avoiding rectilinear minimum Steiner tree.

\begin{lemma} Let $\omega$ be a rectilinearly-convex obstacle that does not have a diagonal. If $\omega$ has a cross $X$, then $X$ is a minimum skeleton for $\omega$.\label{lem-cross}\end{lemma}

\proof $X$ is a connected set of line segments that intersects the extreme edges of $\omega$, and is therefore a skeleton for $\omega$ by Lemma~\ref{lem-connected}. By Lemma~\ref{lem-diagonal}, $\omega$ must have at least two line segments. Therefore $X$ has the smallest possible number of line segments since $|X|=2$.\qed

\begin{lemma}
  Let $\omega$ be an L-obstacle. If $\omega$ has a diagonal $\{d\}$, then $\{d\}$ is a minimum skeleton for $\omega$. Otherwise, a cross is a minimum skeleton for $\omega$.\label{lem-l-obstacle}
\end{lemma}

\proof If $\omega$ has a diagonal $\{d\}$, then $\{d\}$ is a minimum skeleton for $\omega$ by Lemma~\ref{lem-diagonal}. Otherwise, an L-obstacle always has a cross (for instance the two extreme edges whose endpoints are extreme corners) which, by Lemma~\ref{lem-cross} is a minimum skeleton for $\omega$.\qed

\begin{lemma}
  Let $\omega$ be a T-obstacle. Then a minimum skeleton for $\omega$ is any cross of $\omega$.\label{lem-t-obstacle}
\end{lemma}

\proof A T-obstacle does not have a diagonal, but necessarily has a cross $X$ (for example, the pair of edges consisting of the extreme edge $s_1$ of $\omega$ whose endpoints are both extreme corners and a second edge that is an opposite extreme visibility edge orthogonal to $s_1$). Therefore $X$ is a minimum skeleton for $\omega$ by Lemma~\ref{lem-cross}.\qed

Note that minimum skeletons for L-obstacles without diagonals and T-obstacles are not unique; there are infinitely many configurations of crosses that constitute minimum skeletons.

\subsection{Minimum skeletons for staircase obstacles}

While minimum skeletons for rectangles, L-obstacles and T-obstacles have either one or two edges, there is no fixed upper bound (independent of the number of obstacle vertices) on the number of edges in minimum skeletons for the remaining three types in the classification, and constructing minimum skeletons for these types is in general nontrivial. In this section we develop an iterative procedure for constructing minimum skeletons for staircase obstacles. We start by introducing the concept of a maximal extreme visibility edge.

\subsubsection{Maximal extreme visibility edges}

Let $\omega$ be a rectilinearly-convex obstacle and let $S$ be a skeleton for some connected subregion of $\omega$. Then the vertical (horizontal) projection of $S$ is a single vertical (horizontal) line segment. Adding a new edge $s$ to $S$ potentially extends this line segment in one direction or the other, or both. The \emph{vertical (horizontal) advancement} of $s$ is the length of the vertical (horizontal) projection of $s$ that is not covered by the corresponding projection of $S$.

\begin{definition}
  A \emph{maximal extreme visibility edge} $s^*$ is an extreme visibility edge for an extreme edge $e$, such that:
  \begin{itemize}
    \item If $e$ is horizontal, then $s^*$ maximises its vertical advancement and, subject to this, maximises its horizontal advancement towards the other horizontal extreme edge.
    \item If $e$ is vertical, then $s^*$ maximises its horizontal advancement and, subject to this, maximises its vertical advancement towards the other vertical extreme edge.
  \end{itemize}
\end{definition}

\begin{figure}[h]
\centering
\includegraphics[width=14cm]{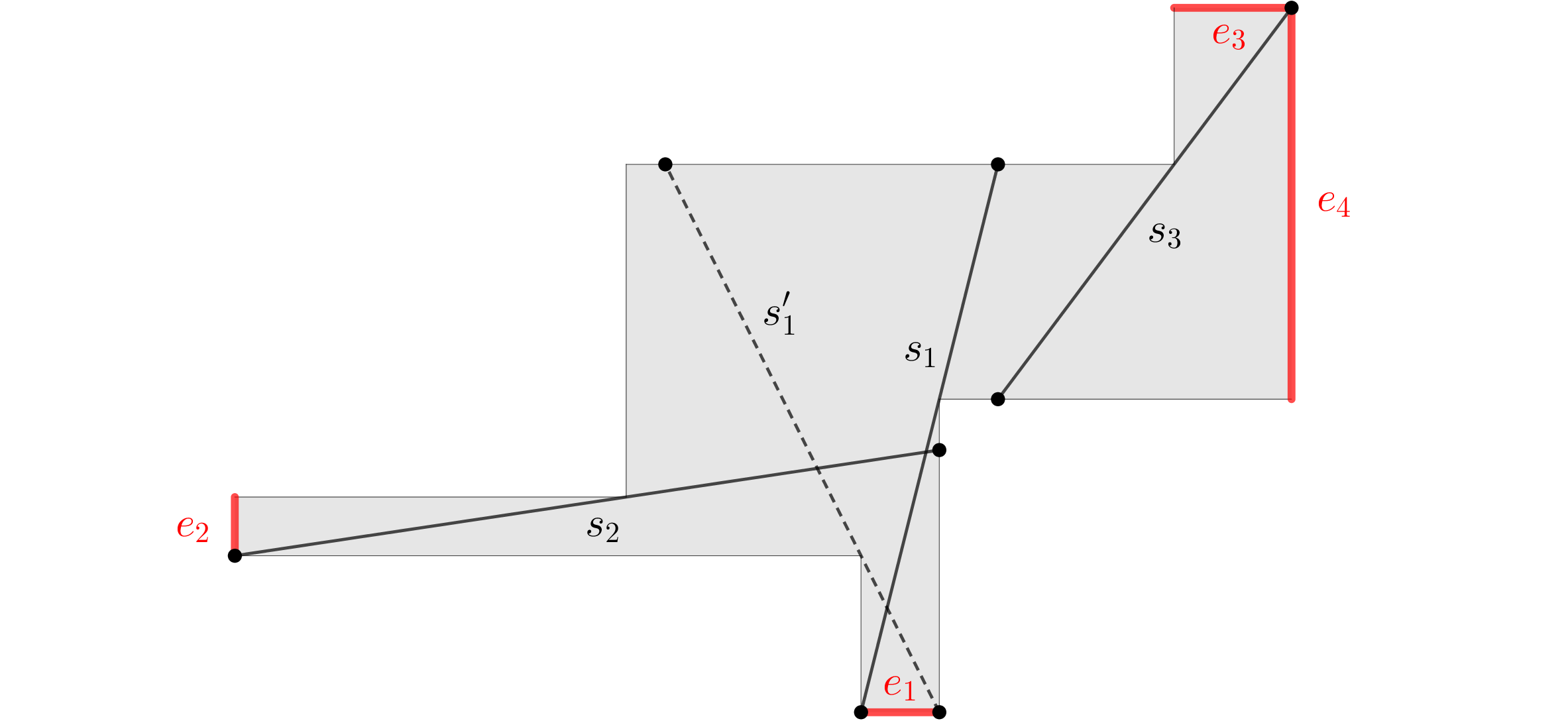}
\caption{A minimum skeleton does not contain $s_1'$ even though the horizontal projection of $s_1'$ is greater than the horizontal projection of $s_1$, since $s_1$ maximises its horizontal advancement towards $e_3$.}\label{fig-advancement}
\end{figure}

The importance of maximising the secondary advancement towards the opposite extreme edge is highlighted in the example of Figure~\ref{fig-advancement} (a partial staircase). It can be shown in this case that $\{s_1, s_2, s_3\}$ is a minimum skeleton for the obstacle shown. If the first edge added to the skeleton is a maximal extreme visibility edge associated with $e_1$, then $s_1$, rather than $s_1'$, should be chosen since $s_1$ has greater horizontal advancement towards $e_3$, even though $s_1'$ has a larger horizontal projection (note that if $s_2$ is added first, then $s_1'$ has no horizontal advancement). Any set of line segments containing $s_1'$ will not be a minimum skeleton.

If a maximal extreme visibility edge has an endpoint at an extreme corner, we refer to it as a \emph{maximal extreme corner visibility edge}. In the following lemma we show that for each extreme corner of a staircase obstacle there exists an associated maximal extreme corner visibility edge.

\begin{lemma} Let $\omega$ be a staircase obstacle with an extreme corner $c$ at the intersection of extreme edges $e_1$ and $e_2$, and let $s_1$ and $s_2$ be maximal extreme visibility edges with endpoints on $e_1$ and $e_2$ respectively. Then $c$ is an endpoint of at least one of $s_1$ and $s_2$.\label{lem-maximal-corner-edge}\end{lemma}

\begin{figure}[h]
\centering
\includegraphics[width=\textwidth]{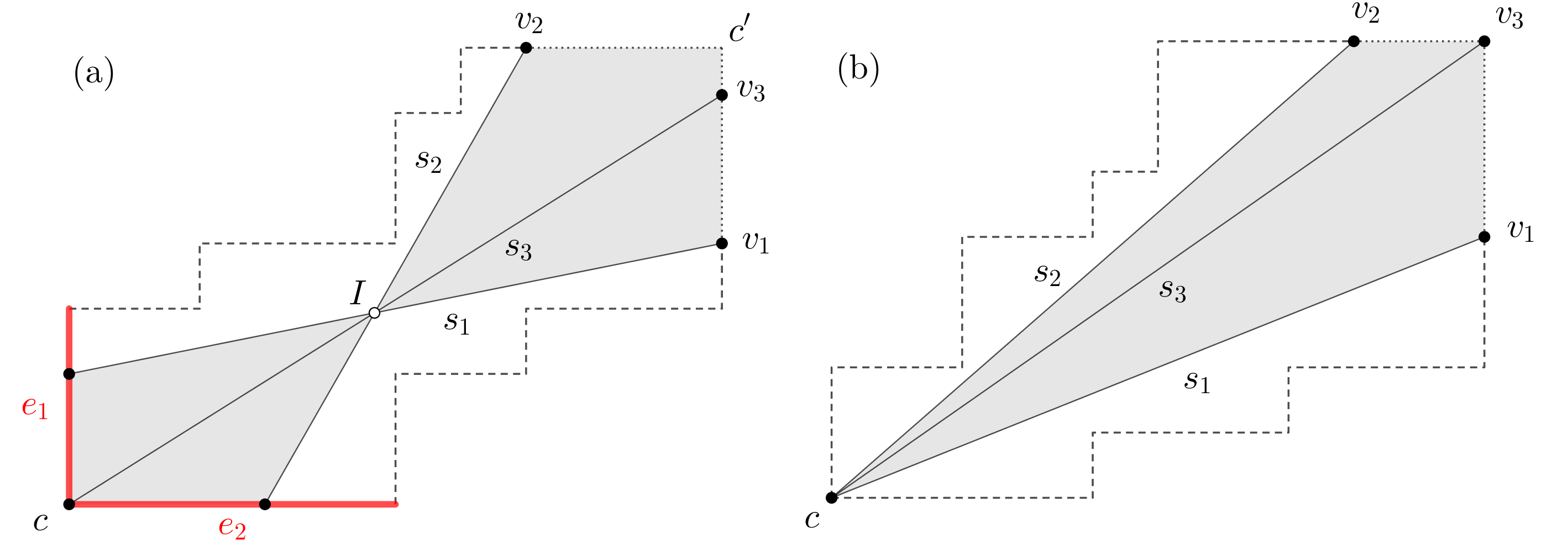}
\caption{Proof of Lemmas~\ref{lem-maximal-corner-edge} and~\ref{lem-corner-maximise}.}\label{fig-extreme-corner-unique-edge}
\end{figure}

\proof We can assume, without loss of generality, that $e_1$ is vertical and $e_2$ is horizontal; let $v_1$ and $v_2$ denote the respective endpoints of $s_1$ and $s_2$ that are opposite $e_1$ and $e_2$. Then (with reference to Figure~\ref{fig-extreme-corner-unique-edge} (a)) $v_1$ must be on the bottom-right staircase walk and $v_2$ on the top-left staircase walk (if $v_1$ is on the top-left staircase walk then it must lie on a horizontal edge of $\omega$, and moving it to the right will increase its horizontal advancement). Therefore the two line segments must intersect at some point $I$. Suppose that neither $s_1$ nor $s_2$ have an endpoint at $c$. Let $c'$ denote the point at the intersection of the vertical line through $v_1$ and the horizontal line through $v_2$. Since $\omega$ is a staircase obstacle, the polygon with vertices $I,v_1,c',v_2$ is inside $\omega$. Let $s_3$ denote the line segment between $c$ and $v_3$, where $v_3$ is the point on the ray from $c$ through $I$ that intersects the rectilinear path $v_1c'v_2$. If $v_3$ is on the line segment $v_1c'$, then $s_3$ has the same horizontal advancement as $s_1$ and has additional vertical advancement, and therefore $s_1$ is not maximal, leading to a contradiction. Otherwise, $v_3$ is on the line segment $v_2c'$, and in this case $s_3$ has the same vertical advancement as $s_2$ and has additional horizontal advancement. Therefore $s_2$ is not maximal, again leading to a contradiction.\qed

In the following lemma we show that a maximal extreme corner visibility edge is unique for a given extreme corner of a staircase obstacle.

\begin{lemma} Let $\omega$ be a staircase obstacle and let $s$ be a maximal extreme corner visibility edge in $\omega$. Then $s$ maximises both its vertical and horizontal advancement.\label{lem-corner-maximise}\end{lemma}

\proof Without loss of generality, assume that $c$ is an extreme corner at the bottom-left corner of $\omega$ (see Figure~\ref{fig-extreme-corner-unique-edge} (b)). Let $s_1$ and $s_2$ be maximal extreme visibility edges which maximise their horizontal and vertical advancements of $\omega$, respectively and, subject to this, maximise their respective vertical and horizontal advancements. Let $v_1$ and $v_2$ denote the endpoints of $s_1$ and $s_2$ that are opposite to $c$, and suppose $v_1$ and $v_2$ are distinct points. Let $v_3$ denote the point at the intersection of the vertical line through $v_1$ and the horizontal line through $v_2$. Since $\omega$ is a staircase obstacle, the polygon with vertices $c,v_1,v_3,v_2$ is contained in $\omega$. Hence the edge $s_3$ between $c$ and $v_3$ lies in $\omega$, and this edge has the same horizontal and vertical advancement as $s_1$ and $s_2$ (respectively), but it has greater advancement in the respective orthogonal directions. Hence $s_1$ and $s_2$ are not maximal, giving a contradiction.\qed

We now show that for any staircase obstacle there exists a minimum skeleton that contains maximal extreme corner visibility edges.

\begin{lemma} Let $\omega$ be a staircase obstacle with extreme corners $c_1$ and $c_2$. Then there exists a minimum skeleton for $\omega$ containing maximal extreme corner visibility edges at $c_1$ and $c_2$.\label{lem-extreme-corner-edge}\end{lemma}

\begin{figure}[h]
\centering
\includegraphics[width=\textwidth]{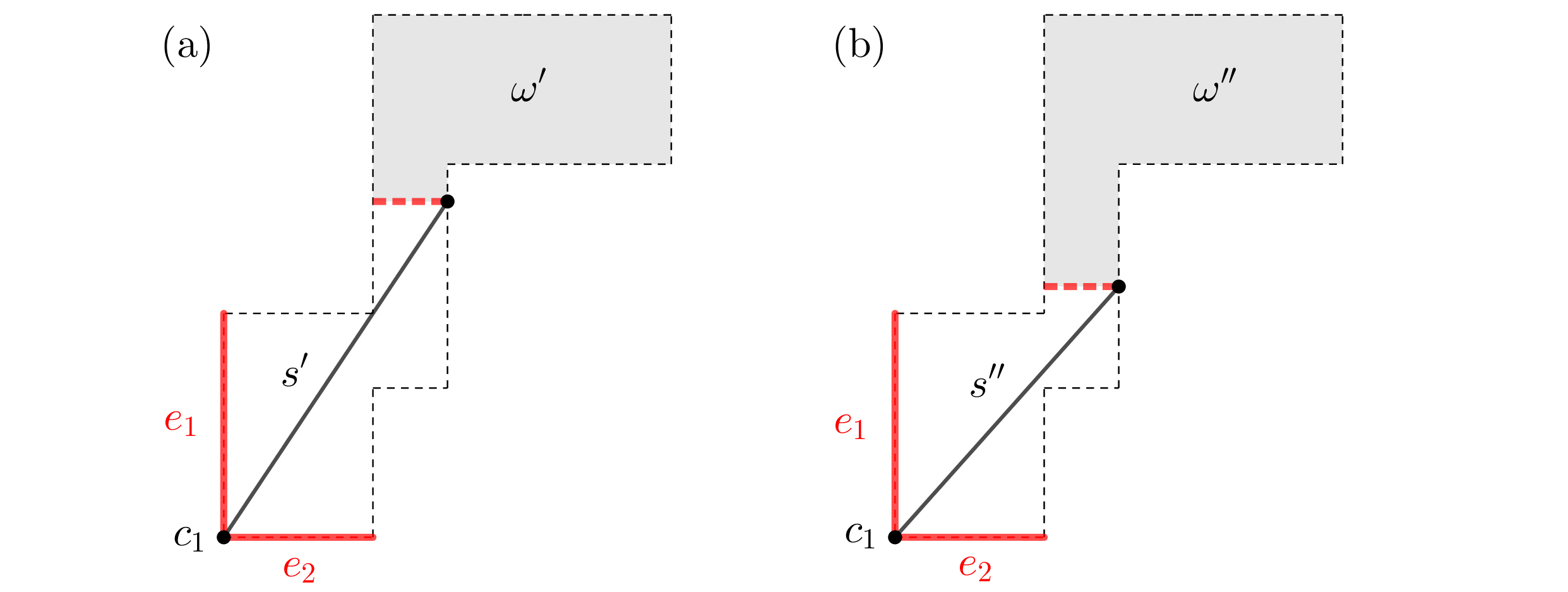}
\caption{Proof of Lemma~\ref{lem-extreme-corner-edge}.}\label{fig-extreme-corner-skeleton-edge-002}
\end{figure}

\proof Let $e_1$ and $e_2$ denote the extreme edges that intersect at $c_1$ (Figure~\ref{fig-extreme-corner-skeleton-edge-002}). By Lemma~\ref{lem-extreme-edge}, both $e_1$ and $e_2$ must be intersected by any minimum skeleton for $\omega$. Let $s_1$ and $s_2$ denote maximal extreme skeleton edges with respective endpoints on $e_1$ and $e_2$. By Lemma~\ref{lem-maximal-corner-edge}, at least one of $s_1$ and $s_2$ has an endpoint at $c_1$. Denote this edge by $s'$. By Lemma~\ref{lem-corner-maximise}, $s'$ is the unique edge that maximises both its horizontal and vertical advancement. Clearly $s'$ is a skeleton for the part of $\omega$ that intersects the bounding box of $s'$, since $s'$ is connected and intersects all four extreme edges of this part of the obstacle.

Let $s''$ be a skeleton edge with an endpoint at $c_1$ that is not maximal. Let $\omega'$ ($\omega''$) denote the part of $\omega$ that remains when the bounding box of $s'$ ($s''$) is subtracted from $\omega$. Let $S'$ ($S''$) denote the minimum sets of edges required to construct skeletons for $\omega$ given the inclusion of $s'$ ($s''$) in the skeleton respectively. Then $|S'|\leq |S''|$ (since $\omega'\subset \omega''$). Therefore $s'$ requires the minimum possible number of additional edges to construct a skeleton for $\omega$, and hence a minimum skeleton for $\omega$ exists that contains $s'$.\qed

\subsubsection{Frontiers and maximal frontier visibility edges}

\begin{figure}[h]
\centering
\includegraphics[width=13.5cm]{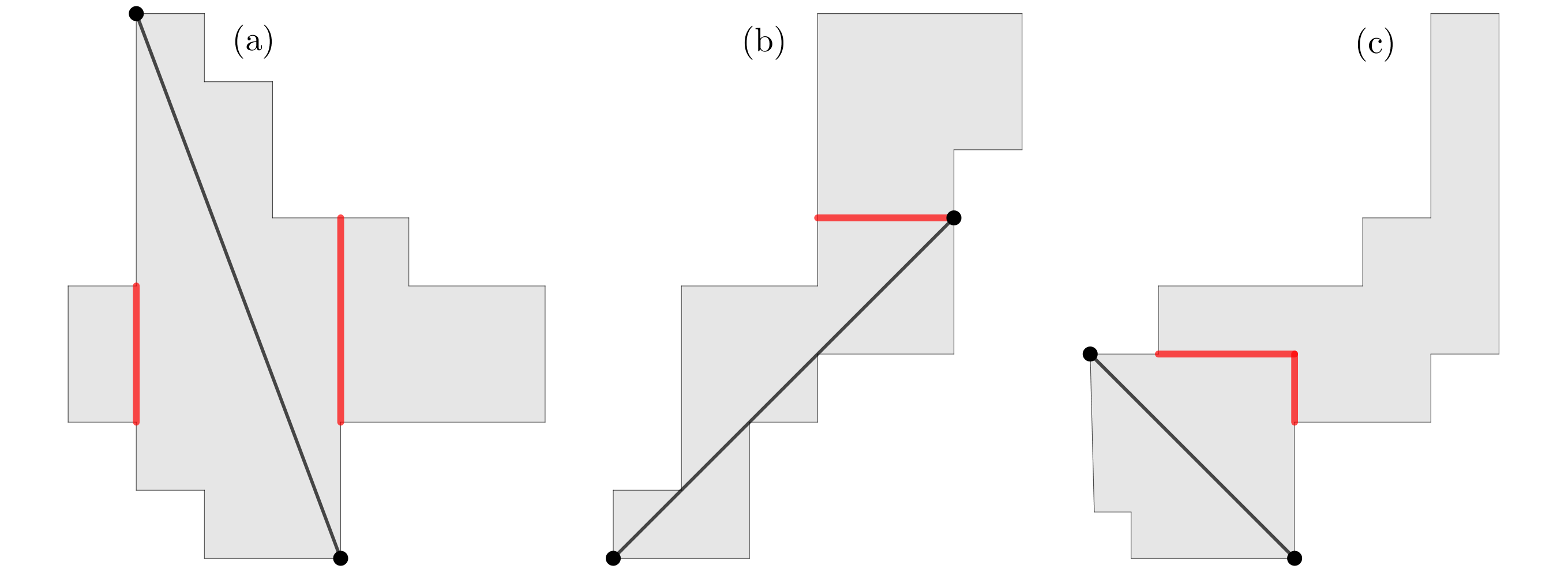}
\caption{Examples of frontiers (frontiers shown as red lines).}\label{fig-examples-of-frontiers}
\end{figure}

Let $S$ be a set of line segments inside a rectilinearly-convex obstacle $\omega$ such that $S$ is a skeleton for $\omega\cap B$, where $B$ is the bounding box of $S$. A \emph{frontier} is the closure of a maximal connected component of the intersection of the boundary of $B$ and the interior of $\omega$. Examples of frontiers are shown in Figure~\ref{fig-examples-of-frontiers}. The function that returns the set of frontiers for a given set of line segments will be denoted by $f(\cdot)$.

A frontier can be a single (horizontal or vertical) line segment, or it can be \emph{L-shaped}, as shown in Figure~\ref{fig-examples-of-frontiers} (c). An L-shaped frontier, which we will treat as a single frontier, can only occur if the intersection of the interior of the frontier with $S$ is empty (provided that all elements of $S$ are maximum length visibility edges). A set $S$ can have multiple associated frontiers (Figure~\ref{fig-examples-of-frontiers} (a)). Frontiers provide a mechanism for constructing skeletons in an iterative fashion. The following lemma provides a useful general property of frontiers which is applicable to any rectilinearly-convex obstacle.

\begin{lemma} Let $\omega$ be a rectilinearly-convex obstacle and let $S^*$ be a minimum skeleton for $\omega$. Let $S\subset S^*$ be a connected or weakly connected set of line segments in $S^*$ and let $F$ be the set of frontiers associated with $S$. Then $f\cap S^*\backslash S\neq \emptyset$ for all $f\in F$.\label{lem-frontiers}\end{lemma}

\proof Suppose that there exists a frontier $f\in F$ that is not intersected by $S^*\backslash S$. Since $\omega$ is convex and given the assumption that all elements of $S$ are maximum length visibility edges, it follows that $f$ divides $\omega$ into two regions, one containing $S$ and the other containing a set $S'\subseteq S^*\backslash S$ that lies on the opposite side of $f$ to $S$. Clearly $S'\cap B(S)=\emptyset$, and therefore $S$ and $S'$ are not weakly connected. This contradicts the assumption that $S^*$ is a skeleton.\qed

A \emph{frontier visibility edge} is a visibility edge that intersects a frontier $f$. A \emph{maximal frontier visibility edge} $s$ is a frontier visibility edge with the following properties:

\begin{itemize}
  \item If $f$ is a vertical line segment, then $s$ maximises its horizontal advancement and, subject to this, maximises its vertical advancement towards the opposite (vertical) extreme edge that is parallel to $f$.
  \item If $f$ is a horizontal line segment, then $s$ maximises its vertical advancement and, subject to this, maximises its horizontal advancement towards the opposite (horizontal) extreme edge that is parallel to $f$.
  \item If $f$ is L-shaped, then $s$ maximises either its horizontal or vertical advancement and, subject to this, maximises its vertical or horizontal advancement (respectively) towards the opposite vertical or horizontal extreme edge.
\end{itemize}

\begin{lemma} Let $\omega$ be a staircase obstacle or a partial staircase obstacle such that $\omega$ has a top-right extreme corner and let $f$ be a frontier associated with a connected or weakly connected set line segments inside $\omega$ that intersects the left and bottom extreme edges of $\omega$. Let $s$ be a maximal frontier visibility edge associated with $f$. Then $s$ maximises both its vertical and horizontal advancement.\label{lem-unique-frontier-edge}\end{lemma}

\begin{figure}[h]
\centering
\includegraphics[width=\textwidth]{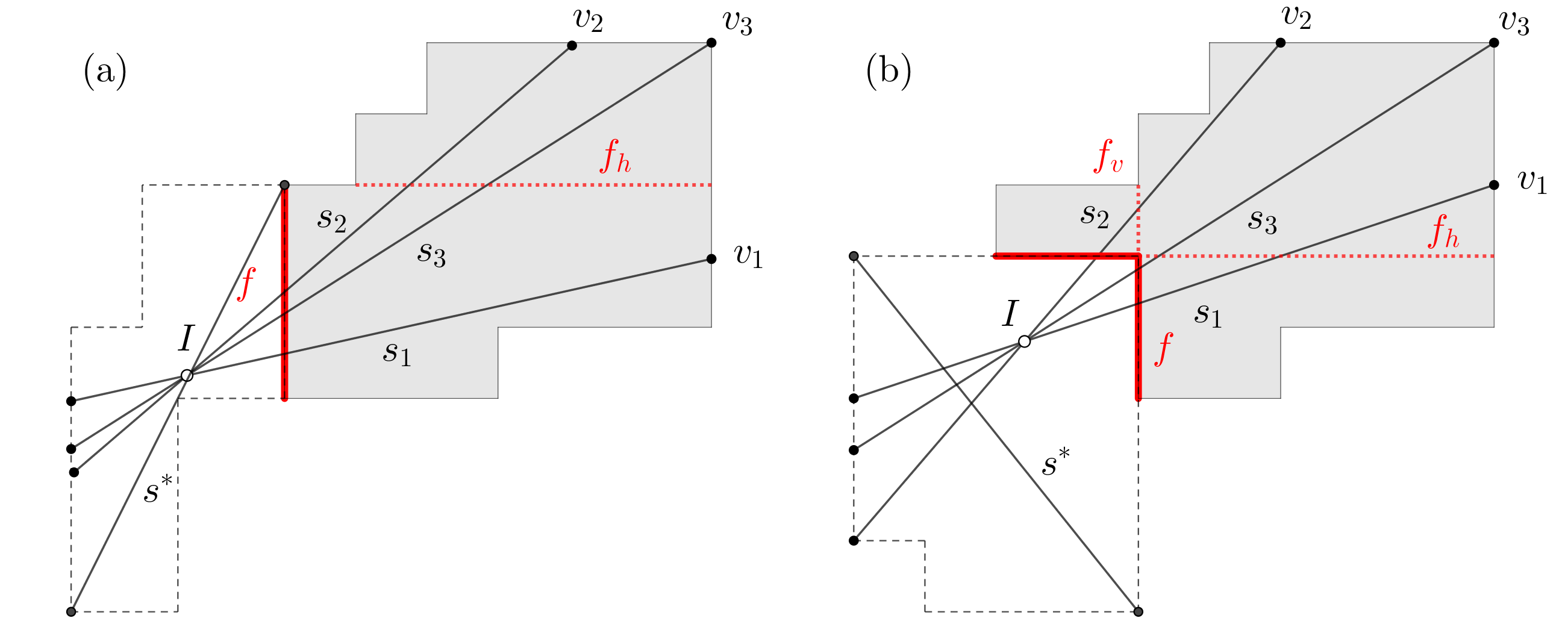}
\caption{Proof of Lemma~\ref{lem-unique-frontier-edge}.}\label{fig-unique-maximal-frontier-edge-2}
\end{figure}

\proof The proof is similar to the proofs of Lemmas~\ref{lem-maximal-corner-edge} and~\ref{lem-corner-maximise}. Assume, without loss of generality, that the direction of advancement is towards the top-right, as in Figure~\ref{fig-unique-maximal-frontier-edge-2}. Let $S$ be a set of skeleton edges for $\omega$ that intersects the bottom and left extreme edges of $\omega$, and let $B$ be the bounding box of $S$. We initially assume that the frontier $f$ corresponding to $B$ is vertical, as in Figure~\ref{fig-unique-maximal-frontier-edge-2} (a). Suppose contrary to the lemma that $s_1$ and $s_2$ are maximal frontier visibility edges which maximise their horizontal and vertical advancements respectively and, subject to this, maximise their respective vertical and horizontal advancements. (Note that in this case the part of $\omega$ below $f_h$, the horizontal line through the topmost point of $f$, has already been covered in the vertical direction.) Let $v_1$ and $v_2$ denote the endpoints of $s_1$ and $s_2$ not in $B$. Let $v_3$ denote the point at the intersection of the vertical line through $v_1$ and the horizontal line through $v_2$. Let $I$ denote the intersection of $s_1$ and $s_2$. Since $s_1$ and $s_2$ are visibility edges (and are therefore in $\omega$), it follows that $I$ must be inside $\omega$, and therefore the polygon with vertices $I,v_1,v_3,v_2$ is inside $\omega$. Let $s_3$ denote the line segment between $v_3$ and the point where the ray from $v_3$ through $I$ intersects the boundary of $\omega$. Hence the edge $s_3$ is inside $\omega$, and this edge has the same horizontal and vertical advancement as $s_1$ and $s_2$ (respectively), and greater advancement in the respective orthogonal directions. Hence $s_1$ and $s_2$ are not maximal, giving a contradiction. Similar arguments apply if the frontier $f$ is horizontal or L-shaped (an example of the latter case is shown in Figure~\ref{fig-unique-maximal-frontier-edge-2} (b)).\qed

\begin{lemma} Let $\omega$ be a staircase obstacle and let $S$ be a set of line segments inside $\omega$ such that $S$ intersects an extreme corner $s_c^*$ of $\omega$ and $S$ is a minimum skeleton for $\omega\cap B(S)$. Let $f$ denote the frontier for $S$, and let $s_f^*$ be the (unique) maximal frontier visibility edge associated with $f$. Then there exists a minimum skeleton $S^*$ for $\omega$ such that $S^*$ contains $S\cup s_f^*$.\label{lem-test}
\end{lemma}

\proof By Lemma~\ref{lem-unique-frontier-edge}, $s_f^*$ maximises both its horizontal and vertical advancement. Therefore $\omega - B(S\cup s_f^*)$ requires the fewest possible number of edges to complete the skeleton. By Theorem~\ref{thm-weakly-connected}, $S\cup s_f^*$ is a skeleton for $B(S\cup s_f^*)$ since $S$ and $s_f^*$ are weakly connected.\qed

\subsubsection{Iterative algorithm for computing minimum skeletons for staircase obstacles}

The preceding results allow us to construct an iterative procedure for computing a minimum skeleton for a staircase obstacle $\omega$ (Algorithm~\ref{algo-compute-staircase}). Note that the algorithm is stated here in a general form in which the inputs include an initial frontier and a termination frontier. This is because the algorithm will be used later as a sub-routine for computing minimum skeletons for partial staircases and general obstacles. For staircase obstacles, the initial frontier and termination frontiers are considered to be points, namely the extreme corners of the obstacles.

\begin{algorithm}
\KwIn{(1) A staircase obstacle $\omega$. (2) An initial frontier $f_0$. (3) A termination frontier $f_t$.}
\KwOut{A minimum skeleton $S^*$ for $\omega$.}
\nl\label{step-test} $S^* := \emptyset$ (initialise the skeleton).\\
\nl $f := f_0$ (initialise the current frontier).\\
\nl $\mathrm{terminate} := \mathrm{FALSE}$ (initialise the termination variable).\\
\nl \While{$\mathrm{\neg terminate}$}{
   \nl $s^* := $ a maximal frontier visibility edge for $f$.\label{algo-staircase-line-5}\\
   \nl $S^* := S^*\cup \{s^*\}$.\\
   \nl \If{$s^*\cap f_t \neq \emptyset$} {
     \nl $\mathrm{terminate} := \mathrm{TRUE}$.
   }
   \nl \Else{\nl $f := $ the new frontier for $S^*$.\label{algo-staircase-line-10}\\
   }
}
\nl \Return{$S^*$}
\caption{\textsc{StaircaseSkeleton}}
\label{algo-compute-staircase}
\end{algorithm}

\begin{figure}[h]
\centering
\includegraphics[width=\textwidth]{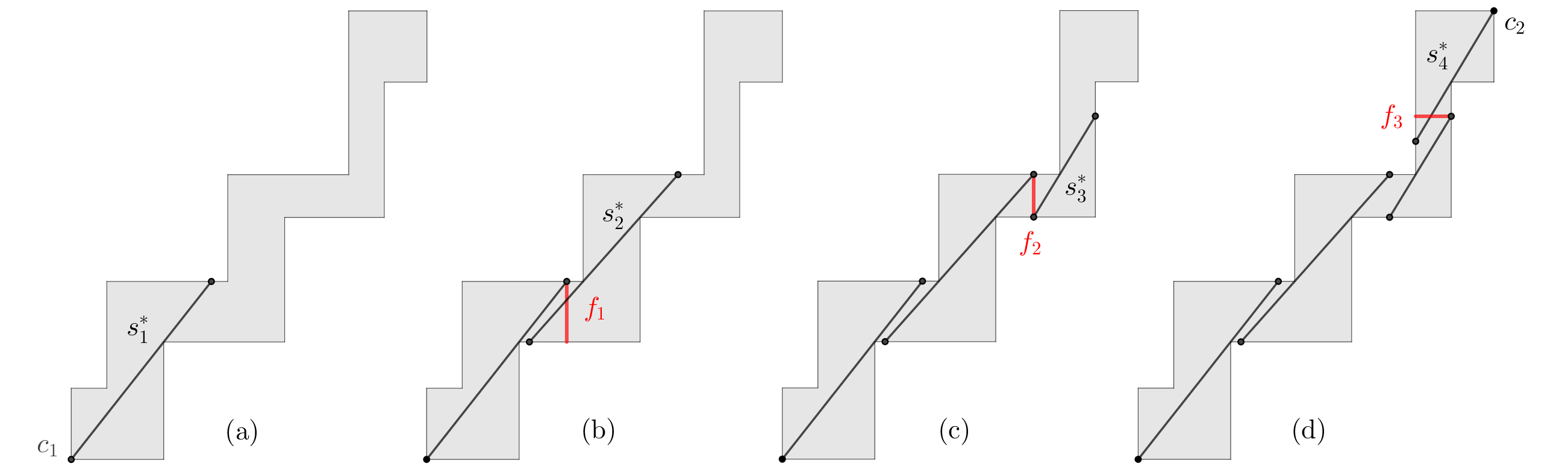}
\caption{Constructing a minimum skeleton using Algorithm~\ref{algo-compute-staircase}.}\label{fig-staircase-examples}
\end{figure}

Methods for efficiently constructing a maximal frontier visibility edge in Line~\ref{algo-staircase-line-5} and the new frontier in Line~\ref{algo-staircase-line-10} are provided in Section~\ref{sect-computational}. Note that Algorithm~\ref{algo-compute-staircase} could alternatively be implemented as a recursive algorithm. We demonstrate the application of Algorithm~\ref{algo-compute-staircase} to the example provided in Figure~\ref{fig-staircase-examples}. At the first iteration (Figure~\ref{fig-staircase-examples} (a)), the (unique) maximal extreme corner visibility edge $s_1^*$ is added to $S^*$ and the frontier $f_1$ associated with $s_1^*$ is computed. In the second and third iterations, the maximal frontier visibility edges $s_2^*$ and then $s_3^*$ are added to $S^*$ and in each case the corresponding frontier is computed. Finally, in the fourth iteration, the maximal frontier visibility edge $s_4^*$ is added to $S^*$. In this case there are a number of candidate frontier skeleton edges that complete the skeleton; our current implementation selects the longest edge from these candidates. The algorithm now terminates since $S^*$ intersects the termination frontier $c_2$.

\begin{theorem} $\textsc{StaircaseSkeleton}(\omega, c_1, c_2)$ computes a minimum skeleton $S^*$ for a given staircase obstacle $\omega$, where $c_1$ and $c_2$ are the extreme corners of $\omega$.\label{thm-staircase}\end{theorem}

\proof The correctness of the algorithm follows from a straightforward inductive argument. By Lemma~\ref{lem-maximal-corner-edge}, a minimum skeleton $S^*$ for $\omega$ exists such that $S^*$ contains a maximal extreme corner skeleton edge $s_1^*$, and by Lemma~\ref{lem-corner-maximise}, $s_1^*$ is unique. The inductive step makes use of Lemma~\ref{lem-test}. At each iteration, the algorithm adds a new maximal frontier visibility edge to the existing set of skeleton edges such that the resulting edge set, $S^*$, is a minimum skeleton for $\omega\cap B(S^*)$.\qed

\subsection{Minimum skeletons for partial staircase obstacles}

\begin{figure}[h]
\centering
\includegraphics[width=12cm]{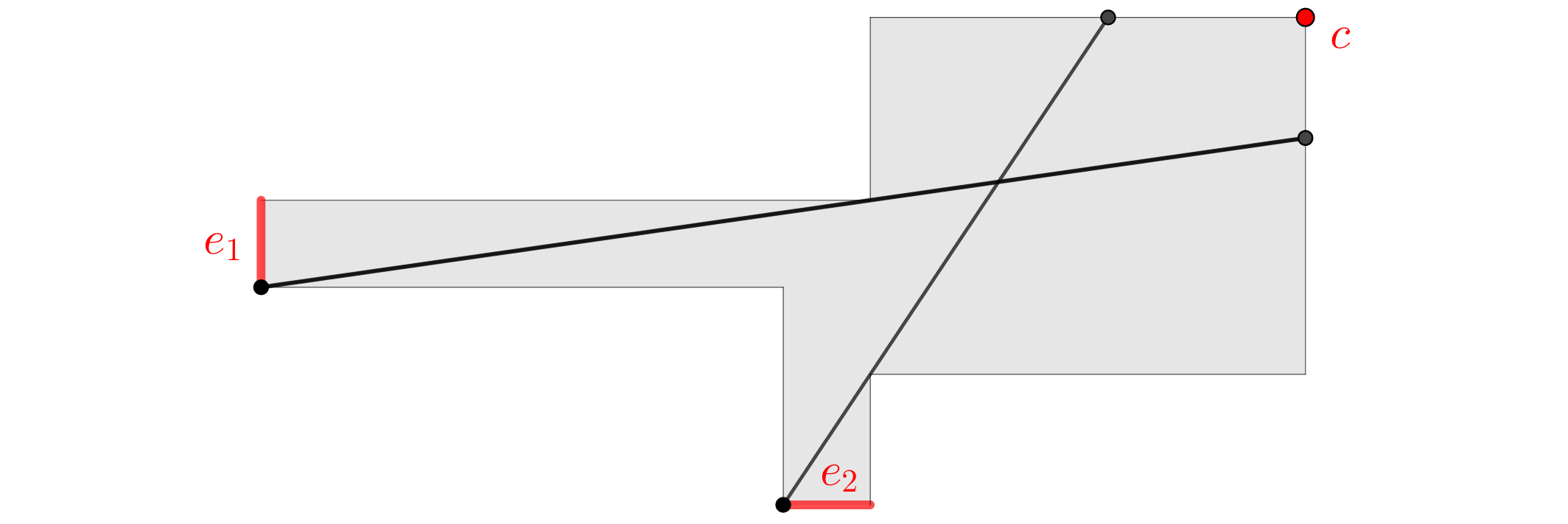}
\caption{A minimum skeleton for a partial staircase $\omega$ does not necessary intersect the extreme corner of $\omega$.}\label{skeleton_shared_extreme_edge}
\end{figure}

A partial staircase has exactly one extreme corner, which will be denoted by $c$. Denote the two extreme edges that do not have an endpoint at $c$ by $e_1$ and $e_2$. By Lemma~\ref{lem-extreme-edge}, both $e_1$ and $e_2$ must be intersected by $S^*$. The example in Figure~\ref{skeleton_shared_extreme_edge} demonstrates that a minimum skeleton for a partial staircase does not necessary intersect $c$. In this case, any edge that has at endpoint at $c$ will be redundant, since $c$ is not visible to $e_1$ or $e_2$, and therefore two additional edges are required with endpoints on $e_1$ and $e_2$ to complete the skeleton (since $e_1$ and $e_2$ are not mutually visible).

The following lemma provides a method for constructing skeleton edges at the end of $\omega$ corresponding to $e_1$ and $e_2$.

\begin{lemma} Let $\omega$ be a partial staircase obstacle with extreme edges $e_1$ and $e_2$ that are opposite to the extreme corner of $\omega$, and assume that $\omega$ does not admit a cross. Let $s_1^*$ and $s_2^*$ be maximal extreme visibility edges with endpoints on $e_1$ and $e_2$ respectively. If $\omega$ has a maximum length adjacent extreme visibility edge $s_{12}^*$ between $e_1$ and $e_2$, then there exists a minimum skeleton $S^*$ such that $s_{12}^*\in S^*$. Otherwise, there exists a minimum skeleton $S^*$ such that $s_1^*\in S^*$ and $s_2^*\in S^*$.\label{lem-partial-ends}
\end{lemma}

\begin{figure}[h]
\centering
\includegraphics[width=\textwidth]{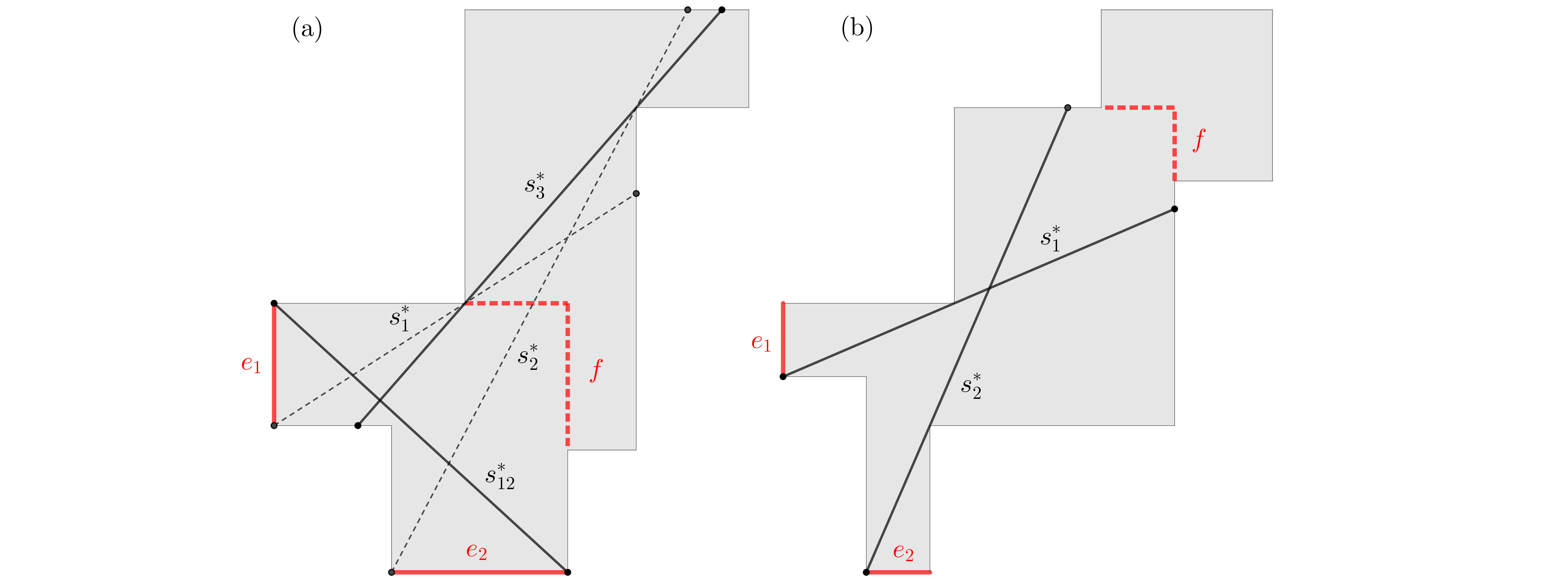}
\caption{Proof of Lemma~\ref{lem-partial-ends}.}\label{figNontrivialEnd}
\end{figure}

\proof Suppose that $\omega$ has a maximum length adjacent extreme visibility edge $s_{12}^*$ between $e_1$ and $e_2$ (Figure~\ref{figNontrivialEnd} (a)). Let $f$ denote the frontier associated with $s_{12}^*$ ($f$ is either L-shaped as in Figure~\ref{figNontrivialEnd} (a), or it is a horizontal or vertical line). Let $s_3^*$ denote the maximal frontier visibility edge for $f$ (where the uniqueness of $s_3^*$ follows from Lemma~\ref{lem-unique-frontier-edge}). Assume, without loss of generality, that the extreme corner of $\omega$ is at the top right of $\omega$, and hence that $s_{12}^*$ is negatively sloped (as in Figure~\ref{figNontrivialEnd} (a)). Then $s_3^*$ can be assumed to be positively sloped, since any line segment intersecting $f$ with a negative slope can be reflected about the vertical line through its midpoint and at least one of the endpoints of the resulting line segment can be extended to increase its horizontal and/or vertical advancement towards $c$. Any positively sloped line segment inside $\omega$ (whose endpoints lie on the boundary of $\omega$) that intersects $f$ must have an endpoint on the bottom-left staircase walk of $\omega$. By the convexity of $\omega$ and since $s_{12}^*$ is a maximum length adjacent extreme visibility edge, we can assume that $s_{12}^*$ is the diagonal of the bounding box for $s_{12}^*\cup e_1\cup e_2$, and therefore $s_{12}^*\cup s_3^*$ is connected. It is also clear that both the horizontal and vertical advancement of $s_3^*$ is at least as great as that of each of $s_1^*$ and $s_2^*$. It follows that $s_{12}^*\cup s_3^*$ is a minimum skeleton for $\omega\cap B(s_{12}^*\cup s_3^*)$, and the region $\omega - B(s_{12}^*\cup s_3^*)$ requires the fewest number of edges to complete the skeleton.

Now suppose that $e_1$ and $e_2$ are not mutually visible (Figure~\ref{figNontrivialEnd} (b)). Using an argument similar to the argument used in Lemma~\ref{lem-unique-frontier-edge}, it can be shown that $s_1^*$ and $s_2^*$ are unique, and since $s_1^*$ and $s_2^*$ are maximal, the region $\omega - B(s_1^*\cup s_2^*)$ requires the fewest number of edges to complete the skeleton.\qed

\begin{figure}[h]
\centering
\includegraphics[width=12cm]{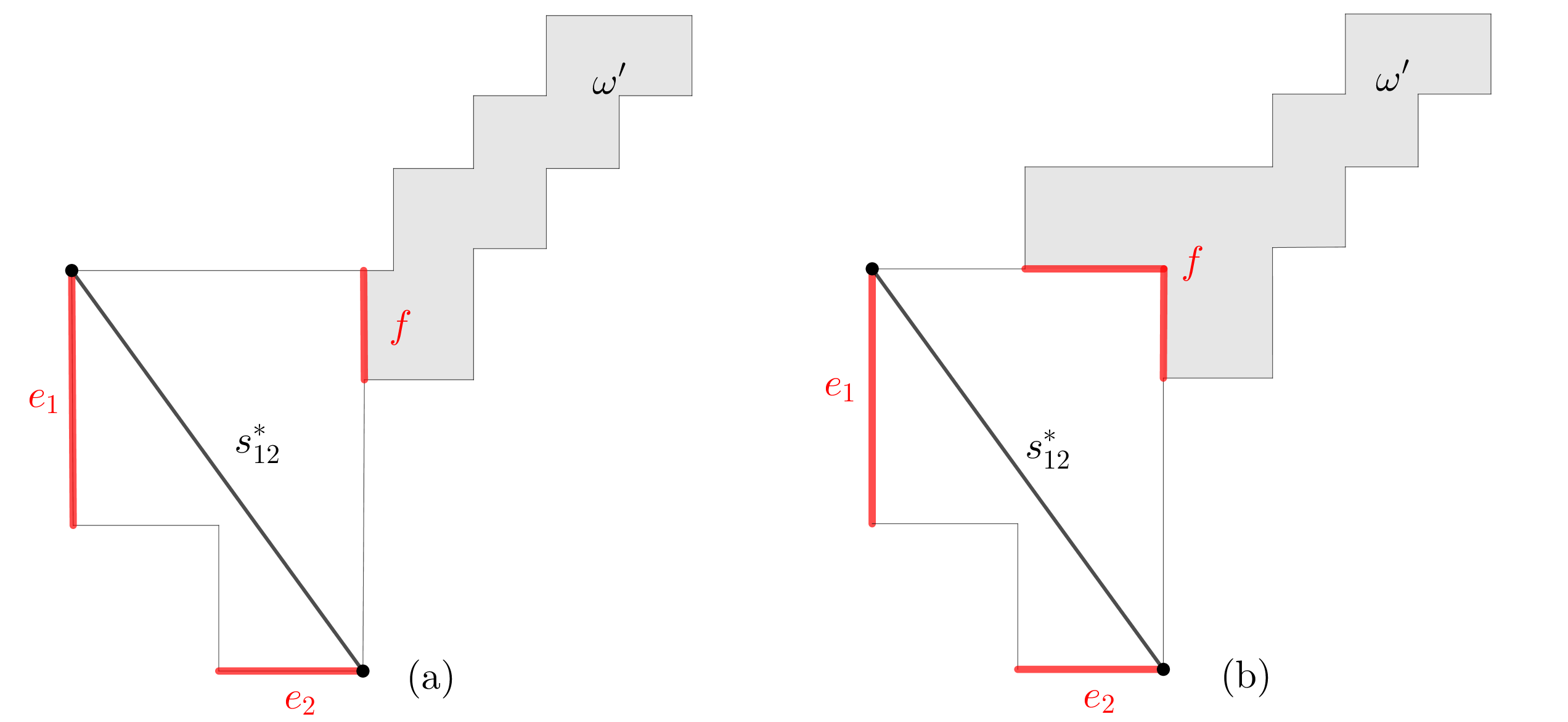}
\caption{(a) $\omega'$ is a sub-staircase with a frontier that is either horizontal or vertical. (b) $\omega'$ can be treated as a sub-staircase with a notch and an L-shaped frontier}\label{figObstacleEnds2}
\end{figure}

Consider a partial staircase obstacle that has a maximum length adjacent extreme visibility edge $s_{12}^*$, and let $B(s_{12}^*)$ denote the bounding box of $s_{12}^*$. When $B(s_{12}^*)$ is subtracted from $\omega$, the remaining obstacle $\omega'$ has two possible forms (see Figure~\ref{figObstacleEnds2}): (a) $\omega'$ is a sub-staircase with a frontier that is either horizontal or vertical. (b) $\omega'$ can be treated as a sub-staircase with a `notch' removed from it and an \emph{L-shaped frontier}. Using the same argument used in the proof of Lemma~\ref{lem-partial-ends}, it can be seen that in both cases a maximal frontier visibility edge associated with $f$ necessarily passes through $f$ and intersects $s_{12}^*$, thereby forming a connected network. We now show that if the minimum skeleton for $\omega$ is not a cross, then there exists a minimum skeleton that intersects the extreme corner of $\omega$.

\begin{lemma} Let $\omega$ be a partial staircase obstacle. If $\omega$ has a cross $X$, then $X$ is a minimum skeleton for $\omega$; otherwise, there exists a minimum skeleton $S^*$ for $\omega$ such that $S^*$ intersects the extreme corner $c$ of $\omega$.\label{lem-partial-staircase-extreme-corner}
\end{lemma}

\proof A partial staircase does not contain a diagonal, since by definition it has exactly one extreme corner. Therefore $|S^*|\geq 2$ by Lemma~\ref{lem-diagonal}. If $\omega$ has a cross $X$, then $X$ is a skeleton for $\omega$ and $|X| = 2$; hence $X$ is a minimum skeleton for $\omega$. Now suppose that $\omega$ does not contain a cross. When the bounding box of the extreme visibility edge or extreme visibility edges that are incident to $e_1$ and $e_2$ (i.e. either $s_{12}^*$ or $s_1^*\cup s_2^*$ using the notation of Lemma~\ref{lem-partial-ends}) is subtracted from $\omega$, then a sub-staircase $\omega'$ is obtained. The argument used in the proof of Lemma~\ref{lem-maximal-corner-edge} can be applied to $\omega'$ to show that there exists a minimum skeleton $S^*$ for $\omega$ such that $S^*$ intersects the extreme corner $c$ of $\omega$.\qed

Algorithm~\ref{algo-compute-partial-staircase} provides a procedure for computing a minimum skeleton for a partial staircase obstacle.

\begin{algorithm}[h]
\KwIn{A partial staircase obstacle $\omega$.}
\KwOut{A minimum skeleton $S^*$ for $\omega$.}
\nl \If{$\omega$ has a cross $X$} {\label{partial-line-1}
    \nl \Return $X$.\label{line-x}
}
\nl \Else{
    \nl Let $e_1,e_2,e_3,e_4$ denote the extreme edges of $\omega$, such that the extreme corner $c$ of $\omega$ is at the intersection of $e_3$ and $e_4$.\\
    \nl \If{$e_1$ and $e_2$ are mutually visible} {\label{partial-line-5}
        \nl Construct the maximum length adjacent extreme visibility edge $s_{12}^*$.\\
        \nl $S^* = s_{12}^*$.\\
    }
    \nl \Else{
        \nl Construct maximal extreme visibility edges $s_1^*$ and $s_2^*$ for $e_1$ and $e_2$ respectively.\\
        \nl $S^* = s_1^*\cup s_2^*$.\\
    }
    \nl $\omega = \omega - B(S^*)$.\\
    \nl $f = f(S^*)$.\\
    \nl \Return $S^*\cup\textsc{StaircaseSkeleton}(\omega, $f$, $c$)$.
}
\caption{\textsc{PartialStaircaseSkeleton}}
\label{algo-compute-partial-staircase}
\end{algorithm}

An efficient procedure for determining if two edges are mutually visible (as required in Lines~\ref{partial-line-1} and~\ref{partial-line-5}) will be discussed in Section~\ref{sect-computational}.

\begin{theorem} $\textsc{PartialStaircaseSkeleton}(\omega)$ computes a minimum skeleton $S^*$ for a given partial staircase obstacle $\omega$.\label{thm-partial-staircase}\end{theorem}

\proof There are three possibilities for $\omega$ that need to be considered (note that by definition a partial staircase obstacle does not have a diagonal):

\begin{enumerate}
  \item $\omega$ has a cross: The existence, or otherwise, of a cross can be determined by checking for the existence of opposite extreme visibility edges. The correctness of Line~\ref{line-x} of the algorithm follows from Lemma~\ref{lem-cross}.
  \item $\omega$ does not have a cross and has a maximum length adjacent extreme visibility edge $s_{12}^*$: By Lemma~\ref{lem-partial-ends} there exists a minimum skeleton $S^*$ such that $s_{12}^*\in S^*$. After constructing the frontier $f$ associated with $s_{12}^*$, the remainder of the skeleton can be constructed using \textsc{StaircaseSkeleton} with initial frontier $f$ and termination frontier $c$.
  \item $\omega$ does not have a cross or a maximum length adjacent extreme visibility edge: By Lemma~\ref{lem-partial-ends}, there exists a minimum skeleton $S^*$ such that $s_1^*\in S^*$ and $s_2^*\in S^*$. After constructing the frontier $f$ associated with $s_1^*\cup s_2^*$, the remainder of the skeleton can be constructed using \textsc{StaircaseSkeleton} with initial frontier $f$ and termination frontier $c$.
\end{enumerate}\qed

\subsection{Minimum skeletons for general obstacles}

Finally, we examine general obstacles (obstacles with no extreme corners), starting with the special cases where $|S^*|=2$ and $|S^*|=3$.

\begin{lemma} Let $\omega$ be a general obstacle and let $S^*$ be a minimum skeleton for $\omega$. Then $|S^*| = 2$ iff $\omega$ has a cross $X$.\label{lem-general-obstacle-cross}\end{lemma}

\begin{figure}[h]
\centering
\includegraphics[width=\textwidth]{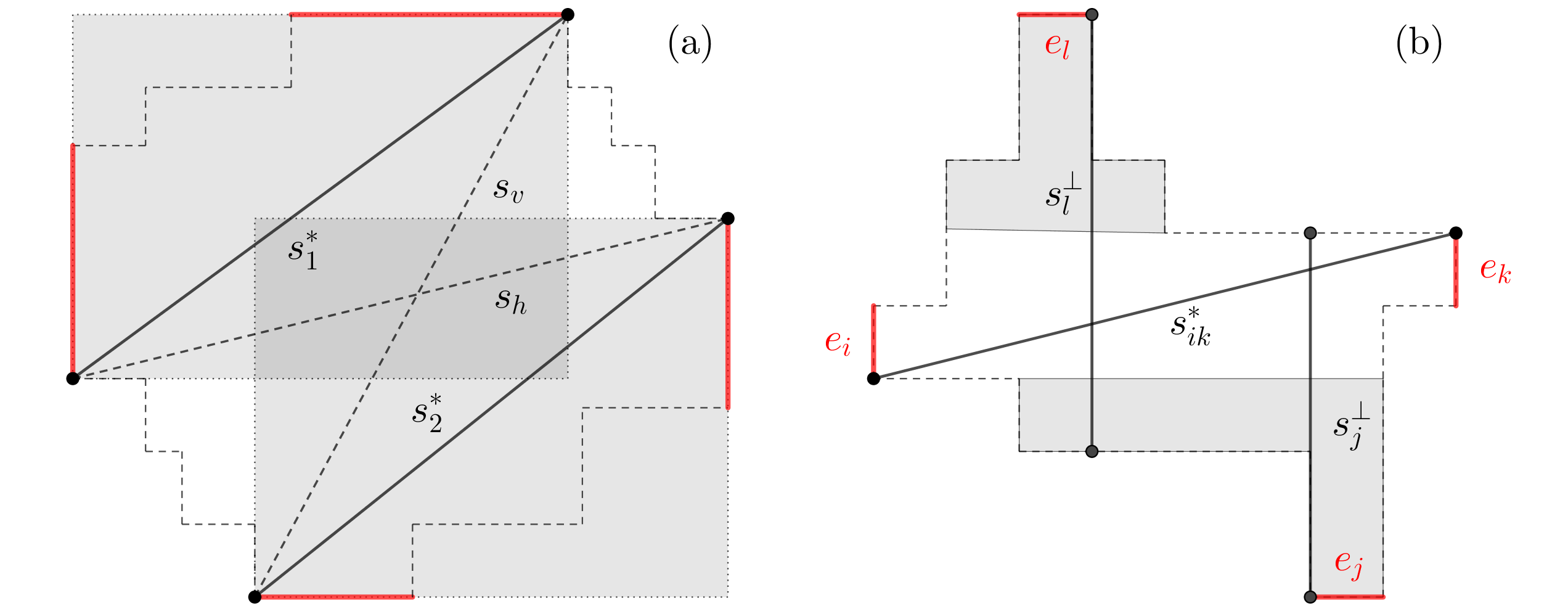}
\caption{(a) An obstacle with a pair of adjacent extreme visibility edges $s_1^*$ and $s_2^*$ also has a cross $s_h$ and $s_v$. (b) A minimum skeleton for an obstacle $\omega$ that has an opposite extreme visibility edge has exactly three edges.}\label{fig-general-obstacle-cross}
\end{figure}

\proof $(\rightarrow)$ If $S^* = 2$, then each extreme edge of $\omega$ contains an endpoint of one of the two edges of $S^*$, say $s_1^*$ and $s_2^*$, by Lemma~\ref{lem-extreme-edge}. Therefore $s_1^*$ and $s_2^*$ are either both opposite extreme visibility edges or they are both maximum length adjacent extreme visibility edges. In the former case $s_1^*$ and $s_2^*$ form a cross. In the latter case, recall that a necessary condition for $S^*$ to be a skeleton is that the projections of $s_1^*$ and $s_2^*$ cover the projections of $\omega$. It follows that the bounding boxes of $s_1^*$ and $s_2^*$ must intersect, and a cross $s_h\cup s_v$ can be constructed between the endpoints of $s_1^*$ and $s_2^*$ (see Figure~\ref{fig-general-obstacle-cross} (a)). $(\leftarrow)$ A cross $X$ is a connected set of edges inside $\omega$ that intersects the four extreme edges of $\omega$. Therefore $X$ is a skeleton for $\omega$ by Lemma~\ref{lem-connected}. It is also a minimum skeleton since $\omega$ has no extreme corners and hence does not have a diagonal, and therefore $|S^*|>1$.\qed

\begin{corollary} Let $S^*$ be a minimum skeleton for a general obstacle that does not admit a cross. Then $|S^*|\geq 3$.\label{lem-triple}\end{corollary}

The following lemma deals with the case when $\omega$ has exactly one opposite extreme visibility edge.

\begin{lemma} Let $\omega$ be a general obstacle that does not admit a cross. If $\omega$ has an opposite extreme visibility edge $s_{ik}^*$ between extreme edges $e_i$ and $e_k$, then $s_{ik}^*\cup s_j^\perp\cup s_l^\perp$ is a minimum skeleton for $\omega$, where $s_j^\perp$ and $s_l^\perp$ are perpendicular extreme visibility edges for $e_j$ and $e_l$ respectively.\label{cor-triple}\end{lemma}

\proof Let $e_j$ and $e_l$ be the two extreme edges of $\omega$ corresponding to $s_j^*$ and $s_l^*$ respectively (Figure~\ref{fig-general-obstacle-cross} (b)). When the bounding box of $s_{ij}^*$ is subtracted from $\omega$, two rectilinearly-convex pieces remain, one containing $e_j$ and the other containing $e_l$. From the convexity of $\omega$, the perpendicular extreme edges $s_j^\perp$ and $s_l^\perp$ clearly intersect $s_{ik}^*$. The three edges together form a connected set of line segments that intersect all four extreme edges of $\omega$ with the least possible number of line segments.\qed

The following lemma addresses the case where $\omega$ has an adjacent extreme visibility edge.

\begin{lemma} Let $\omega$ be a general obstacle with extreme edges $e_i, e_j, e_k, e_l$ such that the ends of $\omega$ (if they exist) are $\{e_i,e_j\}$ and $\{e_k,e_l\}$, and such that $\omega$ does not admit a cross or an opposite extreme visibility edge. If $\omega$ has an adjacent extreme visibility edge $s_{ik}^*$ then $s_{ik}^*\cup s_j^\perp\cup s_l^\perp$ is a minimum skeleton for $\omega$, where $s_j^\perp$ and $s_l^\perp$ are perpendicular extreme visibility edges for $e_j$ and $e_l$ respectively.\label{lem-adjacent-visedge}\end{lemma}

\proof When the bounding box of $s_{ij}^*$ is subtracted from $\omega$, two rectilinearly-convex pieces remain, one piece containing $e_k$ and the other piece containing $e_l$. The perpendicular extreme visibility edges $s_k^\perp$ and $s_l^\perp$ necessarily intersect $s_{ij}^*$ (from the convexity of $\omega$), and the three edges together form a connected set of line segments that intersect all four extreme edges of $\omega$ with the least possible number of line segments.\qed

\subsubsection{Type (c) general obstacles}

The following result is required to construct minimum skeletons for Type (c) general obstacles.

\begin{lemma} Let $\omega$ be a Type (c) general obstacle such that no pair of extreme edges of $\omega$ are mutually visible, and let $s_i^\perp$ denote a perpendicular extreme visibility edge for extreme edge $e_i$. Then $s_1^\perp\cup s_2^\perp\cup s_3^\perp\cup s_4^\perp$ is a minimum skeleton for $\omega$.\label{lem-type-c}\end{lemma}

\proof The proof is similar to the proof of Lemma~\ref{cor-triple}. Suppose that $e_1$, $e_2$, $e_3$ and $e_4$ correspond to the left, bottom, right and top extreme edges respectively and assume that $e_1$ is vertically higher than $e_3$ (and therefore $e_2$ is horizontally to the left of $e_4$). From the convexity of $\omega$, $e_1^\perp$ extends from $e_1$ to a point that is at least as far to the right as the right-most endpoint of $e_4$, and $e_3^\perp$ extends from $e_3$ to a point that is at least as far to the left as the left-most endpoint of $e_2$. Therefore the horizontal projection of $s_1^\perp\cup s_3^\perp$ covers the horizontal projection of $\omega$. When the bounding box of $s_1^\perp\cup s_3^\perp$ is subtracted from $\omega$, two rectilinearly-convex pieces remain, one containing $e_2$ and the other containing $e_4$. From the convexity of $\omega$, the perpendicular extreme visibility edges $s_2^\perp$ and $s_4^\perp$ clearly intersect $e_1^\perp\cup e_3^\perp$. The four edges together form a connected set of line segments that intersect all four extreme edges of $\omega$ with the least possible number of line segments.\qed

\subsubsection{Type (d) general obstacles}

The following lemma provides a result for constructing skeleton edges at the ends of Type (d) general obstacles.

\begin{lemma}
  Let $\omega$ be a general obstacle with ends $\{e_1,e_2\}$ and $\{e_3,e_4\}$ such that $\omega$ does not have a cross, an opposite extreme visibility edge or an adjacent extreme visibility edge between extreme edges at opposite ends of $\omega$. If $\omega$ has an adjacent extreme visibility edge $s_{12}^*$ for the end $\{e_1,e_2\}$, then there exists a minimum skeleton $S^*$ such that $s_{12}^*\in S^*$. Otherwise, there exists a minimum skeleton $S^*$ such that $s_1^*\in S^*$ and $s_2^*\in S^*$, where $s_1^*$ and $s_2^*$ are maximal extreme visibility edges with endpoints on $e_1$ and $e_2$ respectively.\label{lem-skeleton-edges-ends}
\end{lemma}

\begin{figure}[h]
\centering
\includegraphics[width=\textwidth]{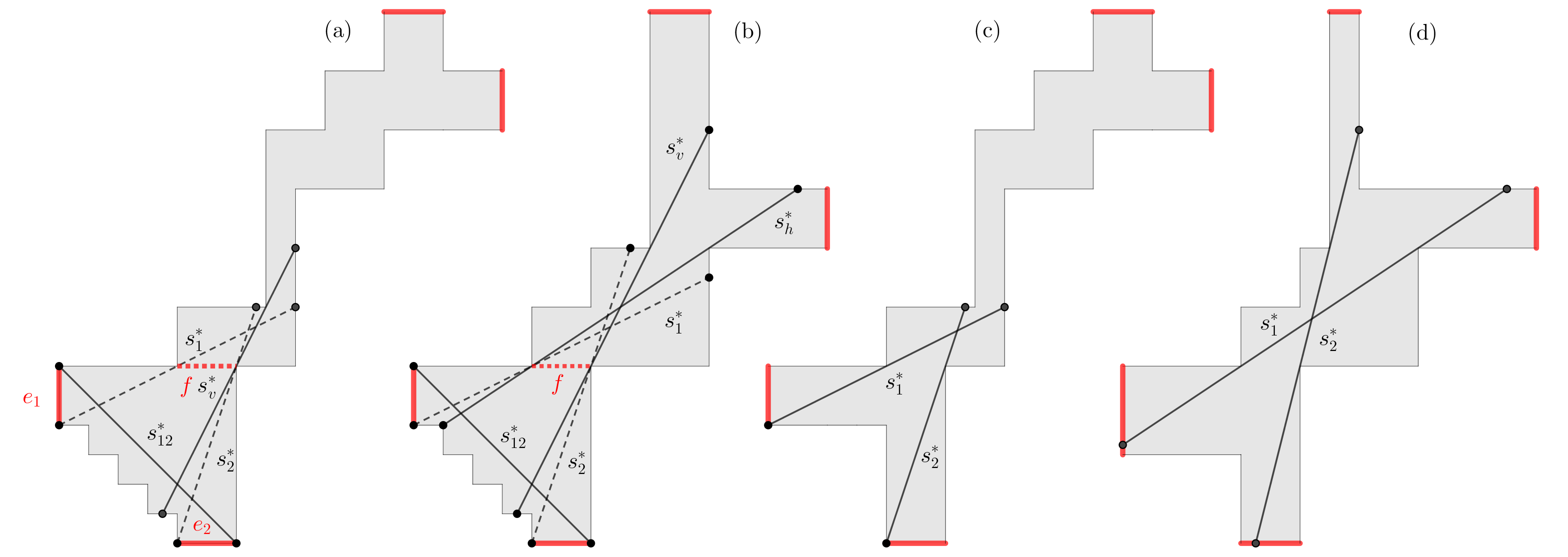}
\caption{Constructing skeleton edges at the ends of a general obstacle. The sub-figure labels correspond to the cases in the proof of Lemma~\ref{lem-skeleton-edges-ends}}\label{fig-general-obstacle-ends}
\end{figure}

\proof Refer to Figure~\ref{fig-general-obstacle-ends}. Without loss of generality, let $e_1$ and $e_2$ denote the left and bottom extreme edges of $\omega$ respectively. Suppose initially that $e_1$ and $e_2$ are mutually visible and let $s_{12}^*$ denote an adjacent extreme visibility edge of $\omega$. Let $f$ denote the frontier associated with $s_{12}^*$. Let $s_1^*,s_2^*,s_3^*,s_4^*$ denote the maximal extreme visibility edges associated with $e_1,e_2,e_3,e_4$ respectively. Let $s_v^*$ and $s_h^*$ denote the two candidates for the maximal frontier visibility edge associated with $f$ (depending on the orientation of $f$), where $s_v^*$ is the line segment intersecting $f$ that maximises its vertical advancement and subject to this then maximises its horizontal advancement, while $s_h^*$ is the line segment intersecting $f$ that maximises its horizontal advancement and subject to this then maximises its vertical advancement. Now we have two cases (illustrated in Figure~\ref{fig-general-obstacle-ends}):

\begin{itemize}
  \item Case (a): Suppose that $s_v^* = s_h^*$ (this is equivalent to saying that the endpoint of $s_v^*$ opposite $e_1$ and $e_2$ does not meet the top-right staircase walk of $\omega$). Then $s_{12}^*\cup s_v^*$ has at least as much horizontal and vertical advancement as $s_1^*\cup s_2^*$, and therefore $\omega - B(s_{12}^*\cup s_v^*)$ requires the fewest possible number of edges to complete the skeleton.
  \item Case (b): Now suppose that $s_v^*\neq s_h^*$ and assume initially that $s_v^*$ has an endpoint on the top-right staircase walk of $\omega$. If $e_3$ and $e_4$ are mutually visible (and have a corresponding adjacent extreme visibility edge $s_{34}^*$), then $s_v^*$ necessarily intersects $s_{34}^*$, and therefore $s_{12}^*\cup s_v^*\cup s_{34}^*$ is a minimum skeleton for $\omega$. If $e_3$ and $e_4$ are not mutually visible, then $s_v^*$ necessarily intersects both of the extreme visibility edges $s_3^*$ and $s_4^*$. In this case $s_{12}^*\cup s_v^*\cup s_3^*\cup s_4^*$ is a minimum skeleton for $\omega$. The preceding arguments are also applicable if $s_h^*$ has an endpoint on the top-right staircase walk of $\omega$. If both $s_v^*$ and $s_h^*$ have an endpoint on the top-right staircase walk of $\omega$, then either of the two edges can be selected for inclusion in $S^*$.
\end{itemize}

Now suppose that $e_1$ and $e_2$ are not mutually visible. Again there are two cases (shown in Figure~\ref{fig-general-obstacle-ends}):

\begin{itemize}
  \item Case (c) If $\omega - B(s_1^*\cup s_2^*)$ is a single region $\omega'$, then $\omega'$ requires the fewest possible number of edges to complete the skeleton.
  \item Case (d) Suppose $\omega - B(s_1^*\cup s_2^*)$ is comprised of two disjoint regions. If $s_{34}^*$ exists then $s_1^*\cup s_2^*\cup s_{34}^*$ is a minimum skeleton for $\omega$; otherwise, $s_1^*\cup s_2^*\cup s_3^*\cup s_4^*$ is a minimum skeleton for $\omega$.
\end{itemize}\qed

Algorithm~\ref{GeneralSkeleton} provides a procedure for computing a minimum skeleton for a general obstacle.

\begin{algorithm}
\KwIn{A general obstacle $\omega$.}
\KwOut{A minimum skeleton $S^*$ for $\omega$.}
\nl Find the extreme edges of $\omega$.\\
\nl For each extreme edge $e_i$, let $s_i^*$ be the corresponding maximal extreme visibility edge and let $s_i^{\perp}$ be any perpendicular extreme visibility edge.\\
\nl \If{$\omega$ has a cross $X$} {\label{line-3-1}
    \nl \Return $X$.
    }
\nl \ElseIf{$\omega$ has an opposite extreme visibility edge $s_{ik}^*$} {
    \nl \Return $s_{ik}^*\cup s_j^\perp\cup s_l^\perp$.
    }
\nl \ElseIf{$\omega$ has an adjacent extreme visibility edge $s_{ik}^*$ and $\omega$ is not a Type (d) obstacle} {
    \nl \Return $s_{ik}^*\cup s_j^\perp\cup s_l^\perp$.\label{line-3-6}
    }
\nl \ElseIf{$\omega$ is a Type (c) obstacle} {
    \nl \Return $s_i^\perp\cup s_j^\perp\cup s_k^\perp\cup s_l^\perp$.
    }\label{line-3-6}
\nl \Else{
    \nl Determine the two ends $\{e_1,e_2\}$ and $\{e_3,e_4\}$ of $\omega$ (as per Figure~\ref{fig-general-obstacles} (d)).\\
    \nl \If{There exists an adjacent extreme visibility edge $s_{14}^*$} {
        \nl \Return $s_{14}^*\cup s_2^\perp\cup s_3^\perp$.
        }
    \nl \ElseIf{There exists an adjacent extreme visibility edge $s_{23}^*$} {
        \nl \Return $s_{23}^*\cup s_1^\perp\cup s_4^\perp$.
    }
    \nl \If{There exists an adjacent extreme visibility edge $s_{12}^*$ for $e_1$ and $e_2$} {
        \nl $S_{12}^* = s_{12}^*$.
    }
    \nl \Else {
        \nl $S_{12}^* = s_1^*\cup s_2^*$.
    }
    \nl \If{There exists an adjacent extreme visibility edge $s_{34}^*$ for $e_3$ and $e_4$} {
        \nl $S_{34}^* = s_{34}^*$.
        }
    \nl \Else {
        \nl $S_{34}^* = s_3^*\cup s_4^*$.
    }
    }
    \nl $S^* = S_{12}^*\cup S_{34}^*$.\\
    \nl \If{$S_{12}^*$ and $S_{34}^*$ are weakly connected} {
    \nl \Return{$S^*$}
    }
\nl \Else {
    \nl $\omega = \omega - B(S_{12}^*) - B(S_{34}^*)$.\\
    \nl $S^* = S^*\cup\textsc{StaircaseSkeleton}(\omega, f(S_{ij}^*), f(S_{kl}^*))$.\\
    \nl \Return{$S^*$}
    }
\caption{\textsc{GeneralSkeleton}}
\label{GeneralSkeleton}
\end{algorithm}

\begin{theorem} $\textsc{GeneralSkeleton}(\omega)$ computes a minimum skeleton $S^*$ for a given general obstacle $\omega$.\label{thm-general-obstacle}\end{theorem}

\proof If $\omega$ has a cross, an opposite extreme visibility edge, an adjacent visibility edge or if $\omega$ is a Type (c) obstacle then $S^*$ can be constructed directly as in Lines~\ref{line-3-1}-\ref{line-3-6} of the algorithm, by Lemmas~\ref{lem-general-obstacle-cross}, \ref{cor-triple}, \ref{lem-adjacent-visedge} and \ref{lem-type-c}. Otherwise, there exists a minimum skeleton $S^*$ such that $S_{12}^*$ and $S_{34}^*$ are elements of $S^*$ by Lemma~\ref{lem-skeleton-edges-ends}. If $S_{12}^*$ and $S_{34}^*$ are weakly connected, then $S_{12}^*\cup S_{34}^*$ is a minimum skeleton for $\omega$ since $|S^*| > 2$, by Lemma~\ref{lem-skeleton-edges-ends}. Otherwise, the remainder of the skeleton can be constructed using \textsc{StaircaseSkeleton} with initial frontier $f(S_{12}^*)$ and termination frontier $f(S_{34}^*)$.\qed

\section{Exact algorithm for computing minimum skeletons}

Algorithm~\ref{ComputeSkeleton} provides an exact algorithm for computing minimum skeletons for rectilinearly-convex obstacles.

\begin{algorithm}
  \KwIn{A rectilinearly-convex obstacle $\omega$.}
  \KwOut{A minimum skeleton $S^*$ for $\omega$.}
  \nl Find the extreme edges and extreme corners of $\omega$, and hence classify $\omega$.\label{line-4-2}\\
  \nl \If{$\omega$ is a rectangle} {
    \nl $S^* = \{d\}$, where $d$ is a diagonal of $\omega$.
  }
  \nl \ElseIf{$\omega$ is an L-obstacle} {
    \nl \If{$\omega$ has a diagonal $d$} {
        \nl $S^* = \{d\}$.
    }
    \nl \Else{\nl $S^* = X$, where $X$ is a cross of $\omega$.}
  }
  \nl \ElseIf{$\omega$ is a T-obstacle} {
    \nl $S^* = X$, where $X$ is a cross of $\omega$.
  }
  \nl \ElseIf{$\omega$ is a staircase obstacle} {
    \nl $S^* = \textsc{StaircaseSkeleton}(\omega, c_1, c_2)$, where $c_1$ and $c_2$ are the extreme corners of $\omega$.
  }
  \nl \ElseIf{$\omega$ is a partial staircase obstacle} {
    \nl $S^* = \textsc{PartialStaircaseSkeleton}(\omega)$.
  }
  \nl \Else{
    \nl $S^* = \textsc{GeneralSkeleton}(\omega).$
  }
  \nl \Return{$S^*$}
  \caption{\textsc{ComputeSkeleton}}
  \label{ComputeSkeleton}
\end{algorithm}

Note that the classification noted in Line~\ref{line-4-2} is determined as in Section~\ref{sect-classification}. We can now state the following theorem.

\begin{theorem} Algorithm~\ref{ComputeSkeleton} computes a minimum skeleton $S^*$ for a given rectilinearly-convex obstacle $\omega$.\label{thm-main}\end{theorem}

\proof From the discussion in Section~\ref{sect-classification}, $\omega$ belongs to one of six possible classifications, based on the number and adjacency of extreme corners. Once the classification of the given obstacle has been identified, a minimum skeleton for the obstacle is obtained as follows:

\begin{itemize}
  \item If $\omega$ is a rectangle, then by Corollary~\ref{cor-rectangle} either of the diagonals of the rectangle is a minimum skeleton for $\omega$.
  \item If $\omega$ is an L-obstacle then if $\omega$ has a diagonal $d$ then $S^* = d$; otherwise $S^*$ is a cross by Lemma~\ref{lem-l-obstacle}.
  \item If $\omega$ is a T-obstacle, then any cross of $\omega$ is a minimum skeleton for $\omega$ by Lemma~\ref{lem-t-obstacle}.
  \item The correctness of the final three cases follow by Theorems~\ref{thm-staircase},~\ref{thm-partial-staircase} and~\ref{thm-general-obstacle} respectively.
\end{itemize}\qed

\subsection{Computational discussion}\label{sect-computational}

Although Algorithm~\ref{ComputeSkeleton} is exact and finite, there are numerous computational steps, such as the construction of maximal visibility edges, that need to be implemented in an efficient way that will scale to obstacles with a large number of vertices. In this section we discuss the current implementations of the key computational steps of Algorithm~\ref{ComputeSkeleton}.

\subsubsection{Determining the obstacle type and extreme edges}

The first step of Algorithm~\ref{ComputeSkeleton} is to determine the type of the given obstacle $\omega$, which is provided as an ordered set of vertices on the boundary of $\omega$. This means that during the input process we can immediately identify $x_{\min}$, $x_{\max}$, $y_{\min}$ and $y_{\max}$ and hence the obstacle type and extreme edges.

The obstacle type is determined by computing the corners $(x_\textrm{min}, y_\textrm{min})$, $(x_\textrm{max}, y_\textrm{min})$, $(x_\textrm{max}, y_\textrm{max})$ and $(x_\textrm{min}, y_\textrm{max})$ of the bounding box $B$ of $\omega$, and identifying which (if any) of the corners are vertices of $\omega$. Let $C$ denote the set of the corners of $B$. Then the obstacle type is determined as follows:

\begin{itemize}
  \item If $|\omega\cap C| = 4$, then $\omega$ is a rectangle;
  \item If $|\omega\cap C| = 3$, then $\omega$ is an L-obstacle;
  \item If $|\omega\cap C| = 2$ and the extreme corners of $\omega$ share an $x$ or $y$ coordinate and hence are adjacent, then $\omega$ is an L-obstacle;
  \item Otherwise, if $|\omega\cap C| = 2$ and the extreme corners of $\omega$ are opposite, then $\omega$ is a staircase obstacle;
  \item If $|\omega\cap C| = 1$, then $\omega$ is a partial staircase obstacle;
  \item If $|\omega\cap C| = 0$, $\omega$ is a general obstacle.
\end{itemize}

The left, bottom, right and top extreme edges of $\omega$ are the unique edges in the boundary of $\omega$ whose endpoints both have $x = x_\textrm{min}$, $y = y_\textrm{min}$, $x = x_\textrm{max}$ and $y = y_\textrm{max}$, respectively.

\subsubsection{Pre-computing the candidate set of skeleton edges (excluding auxiliary edges)}\label{sect-vis-edges}

We define an \emph{auxiliary edge} to be a maximal frontier visibility edge with an endpoint that lies on its corresponding frontier. Let $S^*$ be a minimum skeleton that has been computed by Algorithm~\ref{ComputeSkeleton}. Then $S^*$ is composed of the following types of edges: (1) maximal extreme visibility edges (including opposite, adjacent and extreme-corner visibility edges and diagonals); (2) maximal frontier visibility edges that are not auxiliary edges; and (3) auxiliary edges. Let $G$ denote the set of all candidate skeleton edges that are not auxiliary edges. Then $G$ can be constructed as a pre-processing step to Algorithm~\ref{ComputeSkeleton}, while auxiliary edges are constructed in the course of the algorithm.

In order to construct $G$ we require the following property of maximal frontier visibility edges. We will refer to a staircase obstacle as being \emph{positively sloped} if its extreme corners are located at the bottom-left and top-right corners of the obstacle (denoted by $c_1$ and $c_3$ respectively). The top-left and bottom-right staircase walks of the obstacle will be referred to simply as the top and bottom staircase walks respectively.

\begin{lemma}
  Let $\omega$ be a staircase obstacle that is positively sloped and let $e$ be a maximal frontier visibility edge in $\omega$ for some frontier $f$ of $\omega$, such that $e$ is not an auxiliary edge. Then $e$ satisfies the following properties:
  \begin{enumerate}
    \item $e$ passes through two vertices of $\omega$, say $v_i=(x_i,y_i)$ and $v_j=(x_j,y_j)$ such that $x_j > x_i$ and $y_j > y_i$, and $v_i$ and $v_j$ are on different staircase walks; and
    \item If $v_j\neq c_3$ and $v_j'$ is the endpoint of $e$ that is closest to $v_j$, then $v_j'$ is on the same staircase walk as $v_i$.
  \end{enumerate}\label{lem-alternating}
\end{lemma}

\begin{figure}[h]
\centering
\includegraphics[width=14cm]{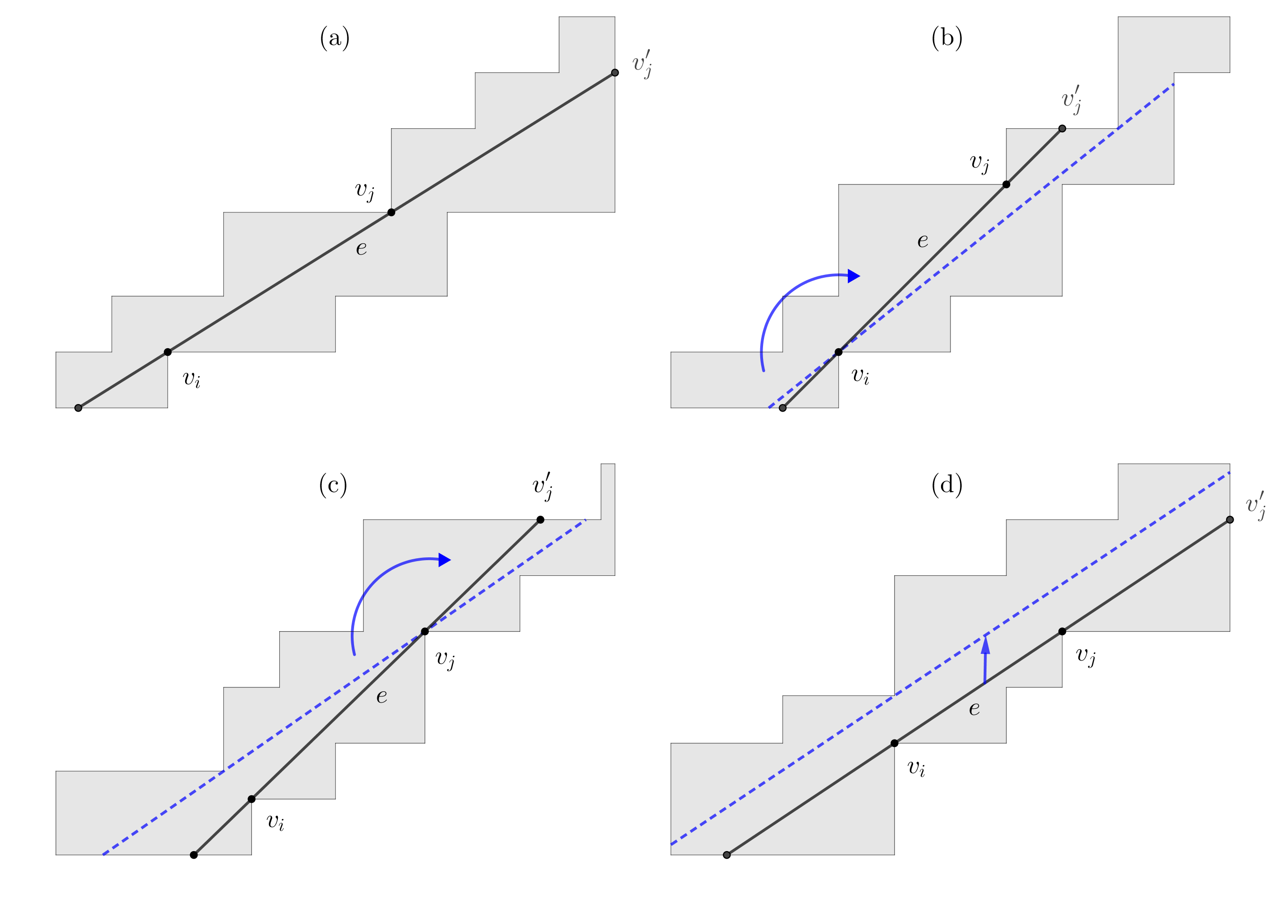}
\caption{(a) A valid skeleton edge has points $v_i$, $v_j$ and $v_j'$ that alternate between opposite staircase walks. (b)-(e) Examples of edges that are not maximal.}\label{fig-alternating-edges}
\end{figure}

\proof Without loss of generality, assume that $v_i$ is on the bottom staircase walk. Then there are four cases, three of which do not satisfy the property of the lemma (Figure~\ref{fig-alternating-edges}):
\begin{itemize}
  \item (a) If $v_j$ and $v_j'$ are on the top and bottom staircase walks respectively, then the conditions of the lemma are satisfied, and there is no continuous transformation of $e$ that increases its advancement.
  \item (b) If $v_j$ and $v_j'$ are both on the top staircase walk then the advancement of $e$ can be increased by rotating $e$ clockwise about $v_i$.
  \item (c) If $v_j$ and $v_j'$ are on the bottom and top staircase walks respectively, then the advancement of $e$ can be increased by rotating $e$ clockwise about $v_j$.
  \item (d) If $v_j$ and $v_j'$ are both on the bottom staircase walk then the advancement of $e$ can be increased by transposing $e$ vertically upwards.
\end{itemize}

In the three cases (a)-(c) in which the candidate edge does not satisfy the property of the lemma, the specified transformation (translation or rotation) increases the horizontal and/or vertical advancement of $e$. Since $e$ intersects the interior of $f$, it will continue to do so under any sufficiently small transformation, and therefore $e$ is not a maximal frontier visibility edge, giving a contradiction.\qed

We now develop an efficient rotational plane sweep method for computing candidate skeleton edges. We begin by looking at maximal frontier visibility edges in staircase obstacles using Lemma~\ref{lem-alternating}.

Let $\omega$ be a staircase obstacle with vertices $V = \{v_1,\ldots,v_n\}$. Assume without loss of generality that $\omega$ is positively-sloped, and that its vertices are labelled in counterclockwise order around its boundary, starting with $v_1$ at the bottom-left extreme corner, and denote the top-right extreme corner by $v_k$. Let $V_b$ and $V_t$ denote the non-convex vertices of $\omega$ on the bottom and top staircase walks (excluding the extreme corners), where a non-convex vertex is a vertex whose interior angle to the obstacle is 270 degrees. The candidate set of skeleton edges for $\omega$ (excluding auxiliary edges and extreme skeleton edges) can be efficiently generated as follows: For each \emph{non-convex vertex} $v_i\in V_b$ do the following (refer to Figure~\ref{fig-rotate-rho}):

\begin{figure}[h]
\centering
\includegraphics[width=13cm]{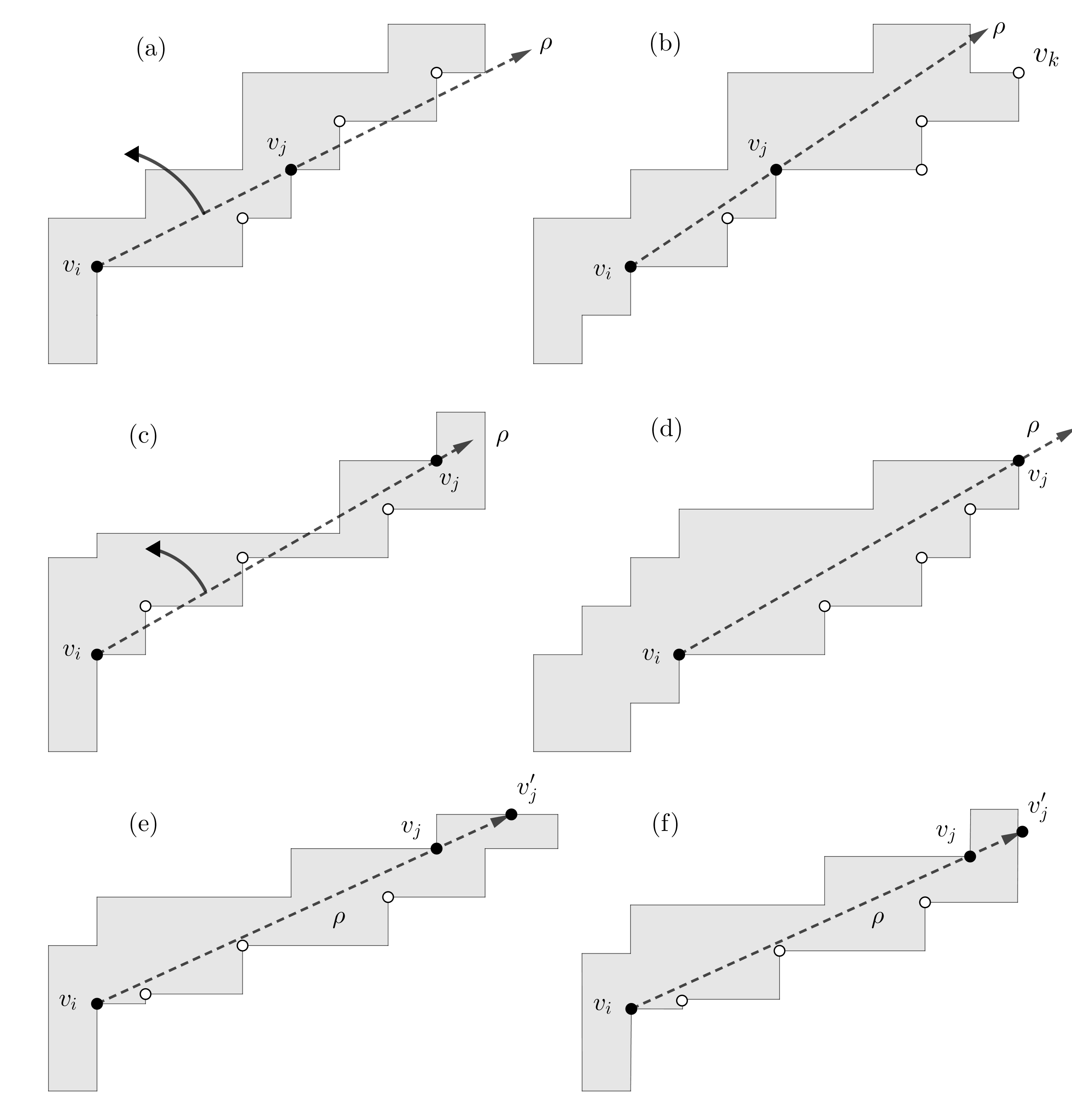}
\caption{Construction of maximal frontier visibility edges.}\label{fig-rotate-rho}
\end{figure}

\begin{enumerate}
  \item Construct a half line $\rho$ starting at $v_i$ that initially points in the positive $x$-direction.
  \item Rotate $\rho$ counterclockwise around $v_i$ until it intersects a non-convex vertex $v_j = (x_j, y_j)$ of $\omega$. There are six possible cases that can occur, each of which is illustrated in Figure~\ref{fig-rotate-rho}:
        \begin{enumerate}
          \item If $v_j\in V_b$ and there exists at least one vertex in $V_b\cup v_k$ that has not already been encountered in the rotational plane sweep by $\rho$, then continue to rotate $\rho$ counterclockwise.
          \item If $v_j\in V_b$ and all other vertices in $V_b\cup v_k$ have already been encountered in the rotational plane sweep by $\rho$, then there are no feasible visibility edges from $v_i$, since the requirements of Lemma~\ref{lem-alternating} cannot be satisfied.
          \item If $v_j\in V_t\cup v_k$ and (1) there exists at least one vertex in $V_b$ that has not been encountered in the rotational plane sweep and whose $y$-coordinate is in $[y_i, y_j]$, or (2) there exists at least one other vertex in $V_t$ that has been previously encountered in the rotational plane sweep and whose $y$-coordinate is in $[y_i, y_j]$, then either continue to rotate $\rho$ counterclockwise, or return no visibility edge if all vertices with $y$-coordinate in $[y_i, y_j]$ have now been encountered in the rotational plane sweep. If continued sweeping is possible, all vertices in $V_t$ whose $y$-coordinate is greater than $y_j$ can be disregarded from the sweep.
          \item If $v_j = v_k$ and all vertices in $V_b$ whose $y$-coordinate is in $[y_i, y_j]$ have already been encountered in the rotational plane sweep, then stop the procedure and do not return a visibility edge (since any visibility edge with an endpoint at $v_k$ will be computed by a separate sweeping procedure applied to the extreme corner).
          \item If $v_j\in V_t$ and all vertices in $V_b$ whose $y$-coordinate is in $[y_i, y_j]$ have already been encountered in the rotational plane sweep, then compute the point $v_j'$ obtained when the line segment between $v_j$ is extended along the line through $v_i$ and $v_j$ to the boundary of $\omega$. If $v_j'$ lies on a horizontal edge of $\omega$ then there are no feasible visibility edges from $v_i$ satisfying the alternating property of Lemma~\ref{lem-alternating}.
          \item Otherwise, $v_j\in V_t$ and all vertices in $V_b$ whose $y$-coordinate is in $[y_i, y_j]$ have already been encountered in the rotational plane sweep, and $v_j'$ lies on a vertical edge of $\omega$, in which case there is a feasible visibility edge that passes through $v_i$ and $v_j$.
        \end{enumerate}
\end{enumerate}

In Case (f), the line segment between $v_i$ and $v_j$ is extended so that its endpoints lie on the boundary of $\omega$. This can be done without intersection computations as follows. To compute $v_j'$, take the set $V_{j'}$ of non-convex vertices of $\omega$ whose $y$-coordinate is greater than $y_j$ and sort the vertices by increasing $y$-coordinate. Search through the vertices by increasing $y$-coordinate until a vertex $v_k$ on the bottom staircase walk is encountered such that $g(v_iv_k) > g(v_iv_j)$, where $g(\cdot)$ denotes the gradient of a line segment. Then $x_{j'} = x_k$ and $y_{j'} = y_i + g(v_iv_j)(x_{j'}-x_i)$. The extension $v_i'$ of $v_i$ is similarly computed, except that non-convex vertices must be checked on the top and bottom staircase walks, since $v_i'$ can be on either side.

A modified version of the above process is also applied to non-convex vertices on the top staircase walk. In this case $\rho$ initially points in the positive $y$-direction, and is rotated in the clockwise direction.

This process is also used to construct other types of skeleton edges, the first three of which can also be constructed during the pre-processing stage:

\begin{itemize}
  \item Maximal extreme corner edges (at the ends of staircase obstacles): The rotational sweep is executed in the clockwise direction with $\rho$ starting at the extreme corner $c_1$. If a feasible edge is not found, then the sweep is executed in the clockwise direction with $\rho$ starting in the positive $y$-direction. To compute the maximal extreme corner visibility edge for $c_3$, the staircase obstacle is reflected about the $x$-axis and the $y$-axis, and the rotational sweep/s performed from $c_3$ (after the reflection $c_3$ is located at the bottom-left corner of $\omega$.
  \item Maximal extreme visibility edges for partial staircase and general obstacles are also computed from one or two rotational plane sweeps (after appropriate reflections have been made to $\omega$).
  \item Adjacent extreme visibility edges are constructed by starting with $v_i$ at an appropriate endpoint of either of the two extreme edges, and performing the rotational plane sweep in the relevant direction.
  \item An auxiliary edge $e_f$ (an edge with an endpoint coinciding with an endpoint $v_f$ of a frontier during the construction of staircase skeletons) is also constructed by applying either a clockwise or a counterclockwise plane sweep with $v_i$ located at $v_f$ (if $e_f$ does not terminate at an extreme corner, then it either passes through a vertex on the top staircase walk and terminates on a vertical edge on the bottom staircase walk, in which case a counterclockwise sweep is required, or it passes through a vertex on the bottom staircase walk and terminates on a horizontal edge on the top staircase walk, in which case a clockwise sweep is required). These edges are constructed during the running of Algorithm~\ref{algo-compute-staircase} as each new frontier is established.
\end{itemize}

\subsubsection{Computing frontiers}

There are two types of frontiers: those associated with maximal frontier visibility edges (or maximal extreme corner visibility edges) computed in the construction of skeletons for staircase obstacles, and those associated with maximum length adjacent extreme visibility edges and pairs of maximal extreme visibility edges (for partial staircase obstacles and general obstacles).

In each case the frontier $f$ is found by computing the closure of the intersection of the bounding box (of the edge or pair of edges) with the interior of $\omega$. Assume without loss of generality that $\omega$ is positively sloped. Now there are three cases:

\begin{enumerate}
  \item If $f$ is for a maximum length adjacent extreme visibility edge, then $f$ is simple to compute in constant time, based on the coordinates of the extreme edges of $\omega$ and their neighbouring edges in $\omega$.
  \item If $f$ is for a maximal frontier visibility edge $s^*$, and assuming that the top-right endpoint $p = (x_p, y_p)$ of $s^*$ is on a vertical edge of $\omega$, then $f$ is a horizontal line segment with one end at $p$ and the other end at $(x_p', y_p)$, where $x_p'$ is the $x$ coordinate of the vertical edge of $\omega$ that lies on the opposite staircase walk to $p$ and whose bottom vertex has the maximum $y$-coordinate among all edges subject to this coordinate being at most $y_p$ (a similar argument applies if $p$ lies on a horizontal edge of $\omega$).
  \item If $f$ is for a pair of maximal extreme visibility edges with respective top-right endpoints $p = (x_p, y_p)$ and $q = (x_q, y_q)$, then (a) if $x_p\geq x_q$ and $y_p\geq y_q$ or $x_q\geq x_p$ and $y_q\geq y_p$, then $f$ is a horizontal or vertical line segment and is computed as in (2) above for $p$ or $q$ respectively; otherwise (b) $f$ is computed as in (1) above based on the corrdinates of the edges of $\omega$ on which $p$ and $q$ lie and their neighbouring edges.
\end{enumerate}

In 2. and 3. (a) it is necessary to sort the list of vertices of $\omega$ by their $x$-coordinates and by their $y$-coordinates.

\subsection{Practical computation of minimum skeletons}

We now show how to efficiently implement the procedure outlined in Algorithm~\ref{ComputeSkeleton} and the other algorithms that it calls. In this section, let $s_i^\perp$ denote any extreme visibility edge for $e_i$ that is perpendicular to $e_i$ and let $s_{ij}^\perp$ denote an opposite extreme visibility edge between $e_i$ and $e_j$ that is perpendicular to $e_i$ and $e_j$.

\subsubsection{Computing minimum skeletons for rectangles, L-obstacles and T-obstacles}

Minimum skeletons for rectangles, L-obstacles and T-obstacles are determined as follows:

\begin{itemize}
  \item If $\omega$ is a rectangle, then either diagonal of $\omega$ is a minimum skeleton for $\omega$.
  \item Let $\omega$ be an L-obstacle with extreme corners $c_1, c_2, c_3$ located at the top-right, bottom-right and bottom-left corners of $\omega$ respectively. If $\omega$ has a diagonal $c_1c_3$, then $c_1c_3$ is a minimum skeleton for $\omega$. Otherwise, the two extreme edges incident with $c_2$ form a minimum skeleton for $\omega$. To determine if $\omega$ has a diagonal, perform a clockwise sweep of an initially-horizontal ray from $c_3$ through $c_2$. If $c_1$ is the first vertex on the top-left staircase walk of $\omega$ encountered, then $c_1c_3$ is a diagonal.
  \item Let $\omega$ be a T-obstacle and let $e_1$ denote the extreme edge of $\omega$ for which both endpoints of $e_1$ are extreme corners of $\omega$. Then $\{e_1, e_{13}^\perp\}$ is a minimum skeleton for $\omega$. Clearly $e_{13}^\perp$ exists for any point on $e_3$ from the convexity of $\omega$.
\end{itemize}

\subsubsection{Computing minimum skeletons for staircase obstacles}

Let $\omega$ be a staircase obstacle and assume without loss of generality that $\omega$ is positively sloped. If $\omega$ has a diagonal $c_1c_3$ then $c_1c_3$ is a minimum skeleton for $\omega$ (the existence of a diagonal can be determined using a similar procedure to the one described above for L-obstacles).

If $\omega$ does not have a diagonal, then the set $G$ of candidate skeleton edges (excluding auxiliary edges) is computed. The extreme corner visibility edge $s_{c_1}^*$ is selected from $G$ and added to the minimum skeleton $S^*$. The frontier $f$ associated with $s_{c_1}^*$ (either a horizontal or vertical line) is computed, and the auxiliary edge $s_a^*$ associated with $f$ is computed and added to $G$. Then the edge in $G$ that intersects $f$ and that has the largest horizontal and vertical advancement from $f$ is selected as the next skeleton edge. The process is repeated until the current frontier intersects the $s_{c_3}^*$, the extreme corner visibility edge at the top-right corner of $\omega$ (at which point $s_{c_3}^*$ is added to $S^*$, completing the skeleton).

\subsubsection{Computing minimum skeletons for partial staircase obstacles}\label{sect-partial}

If necessary, the given partial staircase obstacle is reflected and/or rotated so that its extreme corner $c$ is at the top-right of $\omega$ and $e_1$ and $e_2$ correspond to the left and bottom extreme edges respectively. The set $G$ of candidate skeleton edges (excluding auxiliary edges) is computed and the existence of the adjacent extreme visibility edge $s_{12}^*$ is checked. If $s_{12}^*$ exists, then the (possibly L-shaped) frontier $f$ associated with $s_{12}^*$ is computed, and the edge in $G$ that intersects $f$ with the largest horizontal and vertical advancement from $f$ is added to $S^*$. If $s_{12}^*$ does not exist, then $s_1^*$ and $s_2^*$ are computed, the frontier $f$ associated with $s_1^*\cup s_2^*$ is determined and the maximal edge from $G$ that intersects $f$ is added to $S^*$. The process is repeated until the current frontier intersects the $s_{c_3}^*$, the extreme corner visibility edge at the top-right corner of $\omega$ (at which point $s_{c_3}^*$ is added to $S^*$, completing the skeleton).

\subsubsection{Computing minimum skeletons for general obstacles}

Without loss of generality, let $e_1$, $e_2$, $e_3$ and $e_4$ correspond to the left, bottom, right and top extreme edges, respectively, of a general obstacle $\omega$, and let $x_{\textrm{min}}^i$, $x_{\textrm{max}}^i$, $y_{\textrm{min}}^i$ and $y_{\textrm{max}}^i$ denote the minimum and maximum $x$ and $y$ coordinates of extreme edge $e_i$. As discussed in Section~\ref{sect-general-subtypes} general obstacles can be sub-classified into four sub-types (Figure~\ref{fig-general-obstacles}):

\begin{itemize}
  \item (a) If the projection of $e_1$ onto $e_3$ is not empty and the projection of $e_2$ onto $e_4$ is not empty, then $S^* = s_{13}^\perp\cup s_{24}^\perp$.
  \item (b) Otherwise, if exactly one of the projections (say $e_1$ onto $e_3$) is not empty, then (i) if $e_2$ and $e_4$ are mutually visible, then $S^* = s_{13}^\perp\cup s_{24}^*$; (ii) otherwise, $S^* = s_{13}^\perp\cup s_2^\perp\cup s_4^\perp$.
  \item (c) Otherwise, if $y_{\textrm{min}}^1 > y_{\textrm{max}}^3$ and $x_{\textrm{min}}^4 > x_{\textrm{max}}^2$, then (i) if both pairs of extreme edges are mutually visible, then $S^* = s_{13}^*\cup s_{24}^*$; otherwise (ii) if exactly one pair (say $e_2$ and $e_4$) of extreme edges is mutually visible, then $S^* = s_{24}^*\cup s_1^\perp\cup s_3^\perp$; otherwise (iii) $S^* = s_1^\perp\cup s_2^\perp\cup s_3^\perp\cup s_4^\perp$.
  \item (c) Otherwise, if $y_{\textrm{max}}^1 < y_{\textrm{min}}^3$ and $x_{\textrm{max}}^2 < x_{\textrm{min}}^4$, then (i) if both pairs of extreme edges are mutually visible, then $S^* = s_{13}^*\cup s_{24}^*$; otherwise (ii) if exactly one pair (say $e_2$ and $e_4$) of extreme edges is mutually visible, then $S^* = s_{24}^*\cup s_1^\perp\cup s_3^\perp$; otherwise (iii) $S^*$ is computed by Algorithm~\ref{GeneralSkeleton}, where the obstacle ends are $\{e_1,e_2\}$ and $\{e_3,e_4\}$, and the procedure for efficiently constructing frontiers and maximum frontier visibility edges is essentially the same as that in Section~\ref{sect-partial}.
\end{itemize}

\subsection{Complexity of Algorithm~\ref{ComputeSkeleton}}

The running time of Algorithm~\ref{ComputeSkeleton} is governed by the construction of the set $G$ of candidate visibility edges, where each edge in $G$ is computed by the rotational plane sweep procedure described in Section~\ref{sect-vis-edges}.

Let $\omega$ be a rectilinearly-convex obstacle with $n$ vertices, and let $v$ be a point on the boundary of $\omega$, where $v$ could be a non-convex vertex, an auxiliary point or an endpoint of an extreme edge. The rotational plane sweep procedure requires computing the gradients of the line segments between $v_i$ and every other non-convex vertex (or endpoint of an extreme edge) to the right of $v_i$. The vertices are then sorted by increasing gradient (this sorting can be undertaken in $O(n\log n)$ time using, for instance, the well-known \emph{Heapsort} algorithm~\cite{williams1964}), and a series of constant-time checks are performed on each vertex in the sorted list (see Figure~\ref{fig-rotate-rho}). In the worst case a total of $\frac{n}{2}$ checks are required (i.e. one series of checks for each non-convex vertex of $\omega$).

To construct $G$, the rotational plane sweep procedure is applied to each of the $\frac{n}{2}$ non-convex vertices of $\omega$. In addition, the sweeping procedure is applied for each auxiliary point constructed during the running of Algorithm~\ref{ComputeSkeleton}. We will see in Section~\ref{sect-bounds} that a minimum skeleton has at most $\frac{n}{2}$ edges, and therefore at most $\frac{n}{2}$ auxiliary points. As a consequence, an additional $2\times\frac{n}{2}$ sweeps are required for auxiliary edges (since at most two sweeps are required for each auxiliary edge). From the preceding discussion, we conclude that the overall running time for Algorithm~\ref{ComputeSkeleton} is $O(n^2)$.

Experimental work has demonstrated that the construction of $G$ is fast and scales well in practice as it does not require any intersection tests to be performed to determine the endpoints of candidate skeleton edges.

\subsection{Computing connected skeletons}

\begin{figure}[h]
\centering
\includegraphics[width=\textwidth]{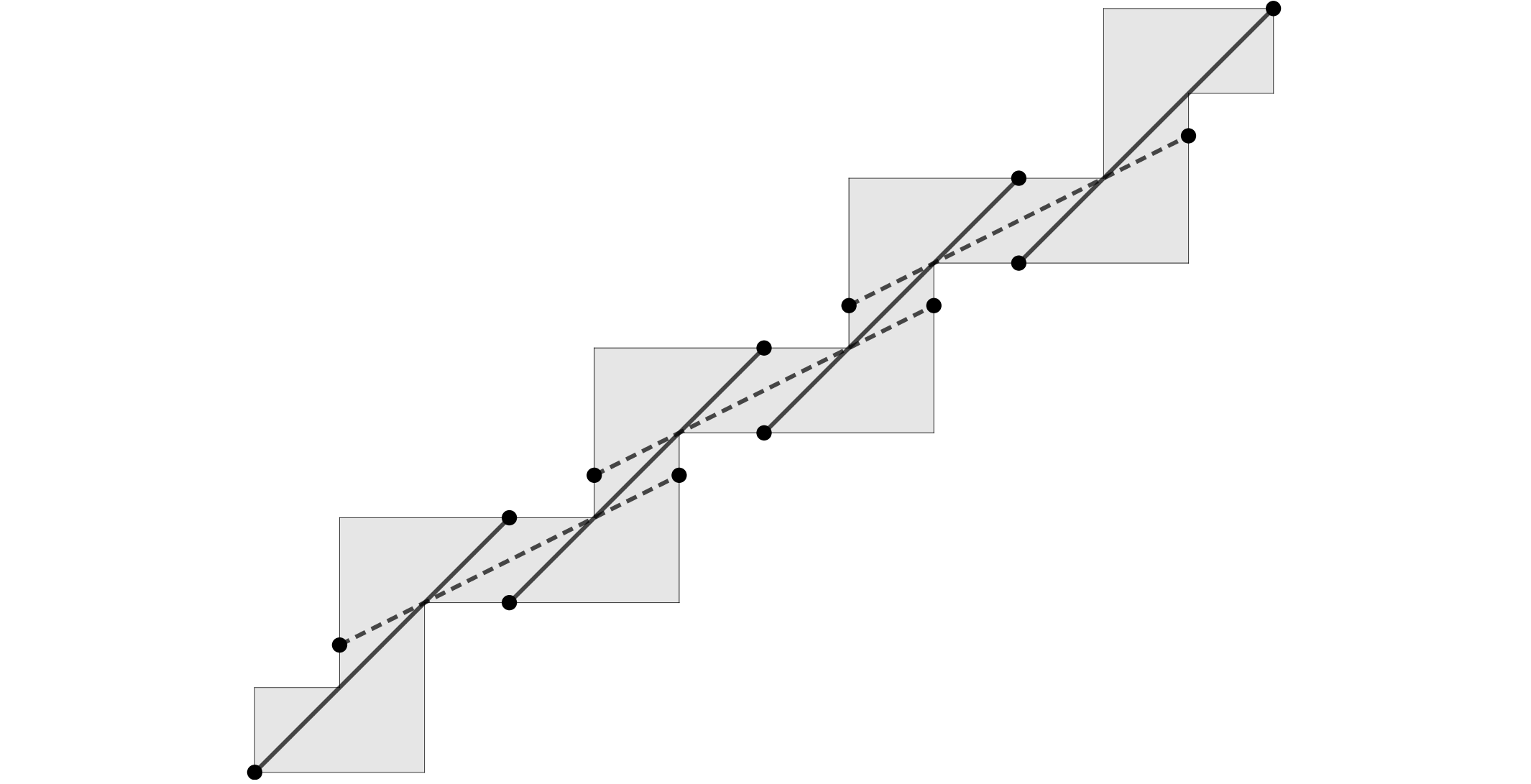}
\caption{Minimum (weakly-connected) skeleton, shown as solid edges. The minimum connected skeleton consists of the solid edges and the dashed edges combined.}\label{fig-connected-bound}
\end{figure}

We define a \emph{connected skeleton} $S$ to be a skeleton for which there is a path between every pair of points in $S$, and a \emph{minimum connected skeleton} to be a connected skeleton with the smallest possible number of line segments. In other words, a connected skeleton is a skeleton for which weak connectivity is not permitted. In this section we show how Algorithm~\ref{ComputeSkeleton} can be modified to compute minimum connected skeletons. In addition to being of theoretical interest, connected skeletons can be faster to compute, because they do not require the added step of generating auxiliary edges during their construction. However, the number of edges in a minimum connected skeleton can be greater than the number of edges in a minimum skeleton that is weakly connected. In fact, the difference between the number of edges in a minimum weakly-connected skeleton and the number of edges in a minimum connected skeleton is unbounded, as demonstrated by the example in Figure~\ref{fig-connected-bound}. In the example, the weakly-connected skeleton $S_1^*$ consists of the four solid edges, while the minimum connected skeleton $S_2^*$ consists of the the same four edges and an additional three edges (shown dashed) to satisfy the connectivity requirement. The structure of the obstacle can be repeated indefinitely such that the resulting obstacle has an arbitrarily large number of vertices, and for any such obstacle we have $|S_2^*| = 2|S_1^*| -1$.

To compute minimum connected skeletons we require the following result.

\begin{lemma}
  Let $\omega$ be a rectilinearly-convex obstacle. Then there exists a minimum connected skeleton $S^*$ for $\omega$ that is a subset of $G$.\label{lem-connected-skeleton}
\end{lemma}

\proof If $\omega$ does not have ends (i.e. if $\omega$ is a rectangle, an L-obstacle, a T-obstacle or a general obstacle of Type (a), (b) or (c), then $S^*$ is constructed from Algorithm~\ref{ComputeSkeleton} without modification, and clearly $S^*$ is a connected subset of $G$. The same is true if $S^*$ has a diagonal, a cross, an opposite extreme visibility edge or an adjacent extreme visibility edge. For obstacles with ends, the skeleton edges at the ends of the obstacle are clearly in $G$. Therefore it is only necessary to look at visibility edges in staircase obstacles (or staircase sub-components of partial staircases and general obstacles that have ends).

\begin{figure}[h]
\centering
\includegraphics[width=\textwidth]{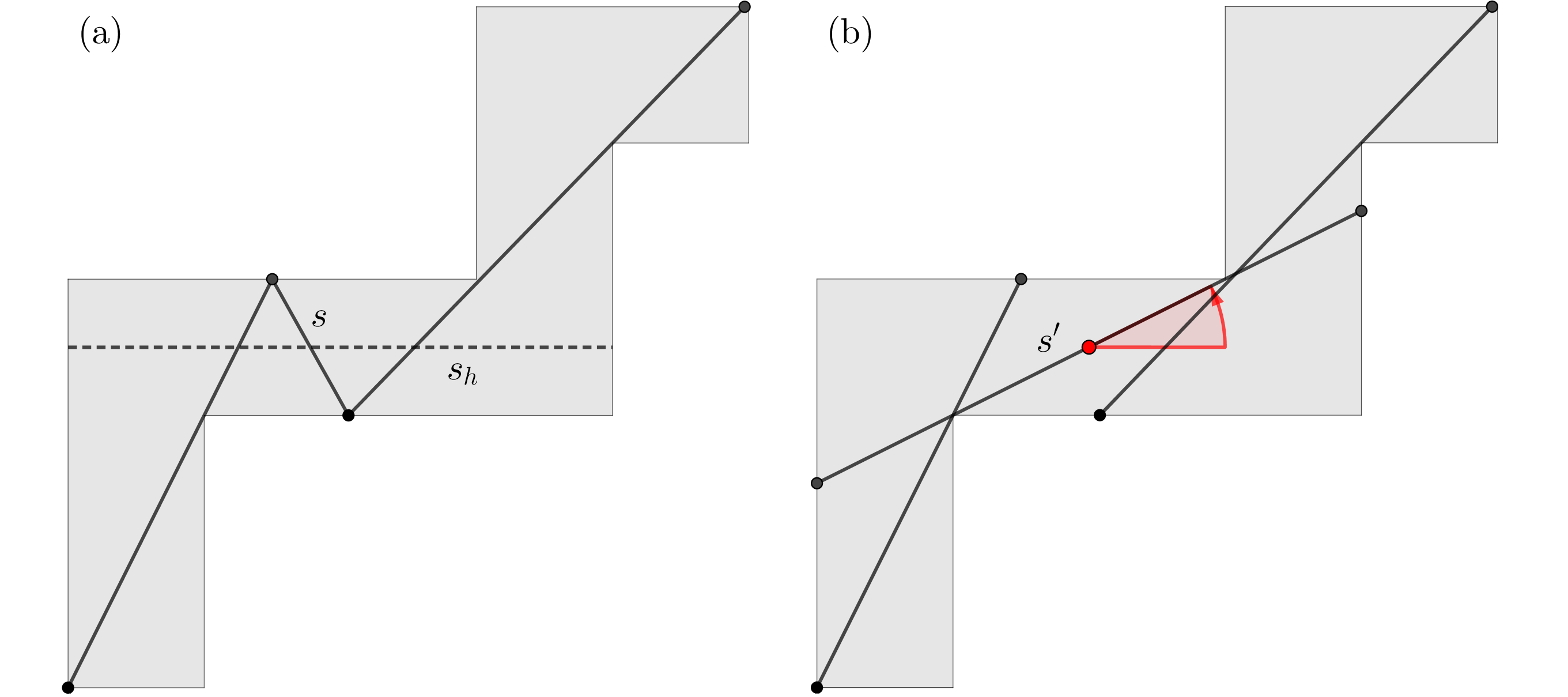}
\caption{Proof of Lemma~\ref{lem-connected-skeleton}.}\label{fig-connected}
\end{figure}

Let $\omega$ be a positively sloped obstacle such that $\omega$ has ends at the bottom-left and top-right corners, and assume that $S^*$ is a minimum connected skeleton for $\omega$. Let $s$ be an edge in $S^*$, and assume that $s$ is negatively sloped (see Figure~\ref{fig-connected} (a)). Then $s$ can be replaced by a horizontal line segment $s_h$ that passes through the midpoint of $s$ and extends as far as possible in both directions to the boundary of $\omega$, and the resulting skeleton remains connected.

Now suppose that $s$ is any edge with a positive (or zero) slope. Then $s$ can be rotated or translated as per Lemma~\ref{lem-alternating} (see Figure~\ref{fig-alternating-edges}) until it is maximal, and all such transformations can be performed without compromising the connectivity of the skeleton. Therefore any connected set of line segments in $S^*$ can be replaced by maximal edges from $G$.\qed

\begin{figure}[h]
\centering
\includegraphics[width=\textwidth]{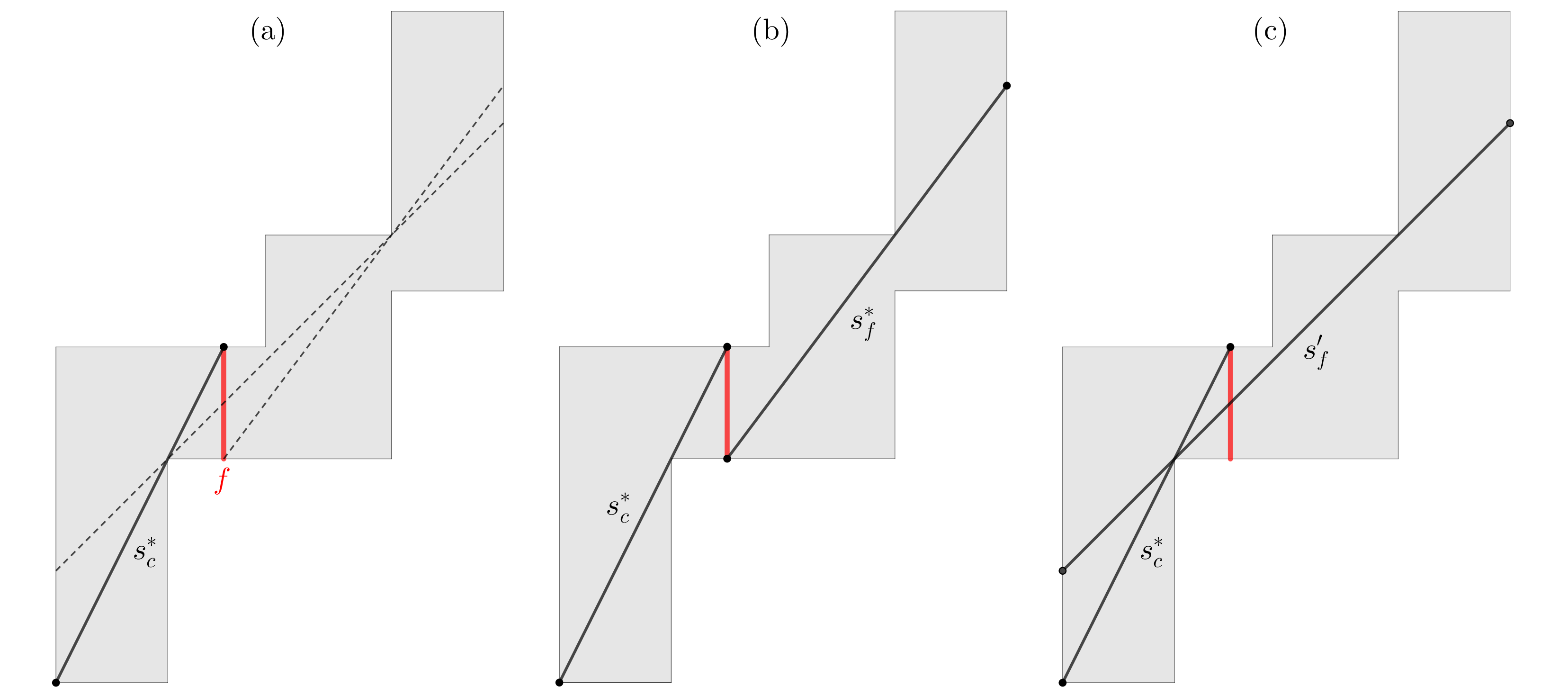}
\caption{An example demonstrating the required modification of Algorithm~\ref{ComputeSkeleton} when the goal is to compute a connected skeleton.}\label{fig-vis-edges-at-frontier}
\end{figure}

A minimum connected skeleton can be constructed by modifying Algorithm~\ref{algo-compute-staircase} so that at each iteration, the maximal frontier visibility edge is an edge from $G$ with the greatest advancement among all edges that intersect the previously constructed skeleton edge $s^*$. The modification is illustrated in Figure~\ref{fig-vis-edges-at-frontier}. At the first iteration (a), the maximal extreme corner visibility edge $s_c^*$ is added to the skeleton. If the goal were to construct a weakly-connected skeleton (b), then the auxiliary edge $s_f^*$ would be added to the skeleton, since it has the greatest advancement among all edges that are weakly connected with $s_c^*$. However, in this sub-section the goal is to construct a connected skeleton, and in this case $s_f'$ is added to the skeleton, since it is has the greatest advancement among all edges that intersect $s_c^*$. The modified algorithm does not require the construction of auxiliary edges, and therefore can potentially have a significantly reduced running time in practice.

For partial staircase obstacles and general obstacles, there are no modifications required when computing the skeleton edges at each end of the obstacle, since any such edges will be connected regardless of whether or not weak connectivity is permitted.

\section{Computational bounds and results}

In this section we show that for any rectilinearly-convex obstacle the number of edges in a minimum skeleton is at most half the number of edges in the boundary of the obstacle, and in most cases is significantly smaller. We have implemented Algorithm~\ref{ComputeSkeleton} and applied it to randomly-generated rectilinearly-convex obstacles.

\subsection{Upper bound on the number of edges in a minimum skeleton}\label{sect-bounds}

Although there is no upper bound (independent of the number of obstacle vertices) on the number of edges in a minimum skeleton for staircase obstacles, partial staircase obstacles and general obstacles, we can nevertheless bound the number of skeleton edges based on the number of vertices of $\omega$ as follows.

\begin{lemma} Let $\omega$ be a rectilinearly-convex obstacle with vertex set $V$ and edge set $E$, and let $S^*$ be a minimum skeleton for $\omega$. Then $|S^*|\leq\frac{|V|}{2} = \frac{|E|}{2}$.\end{lemma}

\proof For any polygonal obstacle it is clear that the number of vertices is the same as the number of edges. If $\omega$ is a rectangle, then $|V|=4$ and $|S^*|=1$. If $\omega$ is an L-obstacle, then $|V|\geq 6$ and $|S^*|\leq 2$. If $\omega$ is a T-obstacle, then $|V|\geq 8$ and $|S^*|=2$.

\begin{figure}[h]
\centering
\includegraphics[width=\textwidth]{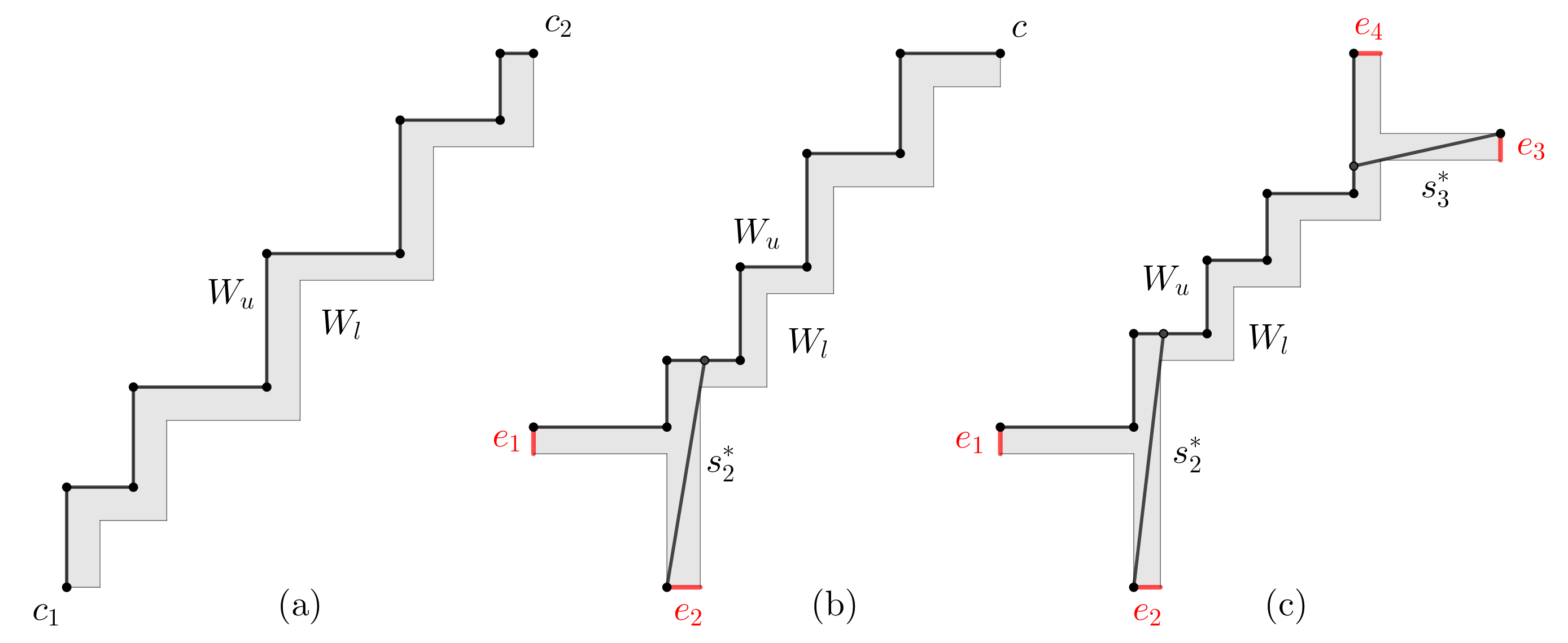}
\caption{Upper bound on the number of edges in a minimum skeleton. (a) Staircase obstacle. (b) Partial staircase obstacle. (c) General obstacle.}\label{fig-max-skeleton-edges}
\end{figure}

Now suppose that $\omega$ is a staircase obstacle and assume without loss of generality that $\omega$ is positively sloped (Figure~\ref{fig-max-skeleton-edges} (a)). Let $W_u$ and $W_l$ denote the sets of edges in the top-left and bottom-right staircase walks of $\omega$ respectively and without loss of generality assume that $|W_u|\leq |W_l|$. Then $W_u$ is a skeleton for $\omega$ since it is a connected set of line segments that intersects the four extreme edges of $\omega$, and $|W_u|$ has at most $\frac{|V|}{2}$ edges since $|W_u| + |W_l| = |E| = |V|$.

Now suppose that $\omega$ is a partial staircase obstacle and without loss of generality assume that $e_1$ and $e_2$ are the (disjoint) left and bottom extreme edges and $c$ is an extreme corner at the top right of $\omega$ (Figure~\ref{fig-max-skeleton-edges} (b)). Let $W_u$ and $W_l$ denote the top-left and bottom-right staircase walks of $\omega$ excluding $e_1$ and $e_2$ respectively, and assume that $|W_u|\leq |W_l|$. Then $W_u\cup s_2^*$ is a skeleton for $\omega$ (where $s_2^*$ is the maximal extreme visibility edge for $e_2$) with $|W_u| + 1$ edges and since $\omega$ has at most $2|W_u| + 4$ edges, we have that $|S|\leq \frac{|V|}{2}-1$. If $e_1$ and $e_2$ are mutually visible, then the same argument applies when $s_2^*$ is replaced by the adjacent extreme visibility edge $s_{12}^*$.

Similar arguments apply to the case where $\omega$ is a general obstacle (Figure~\ref{fig-max-skeleton-edges} (c)). In this case $W_u\cup s_2^*\cup s_3^*$ is a skeleton for $\omega$ with at most $\frac{|V|}{2}-2$ edges.

The bound stated in the lemma is tight and can be achieved by constructing a `skinny' staircase obstacle for which $|W_u|\leq |W_l|$ and $W_u$ and $W_l$ are closely aligned, in which case $W_u$ is a minimum skeleton.\qed

\subsection{Random generation of rectilinearly-convex obstacles}

As far as we are aware, no algorithm exists for generating random rectilinearly-convex polygons. Therefore we have constructed the following procedure for generating random rectilinearly-convex obstacles where the type, dimensions of the bounding box and number of vertices are given (Algorithm~\ref{RandomRectilinearlyConvexObstacle}). Note that rectangles are not included since they are immediately determined by the bounding box.

\begin{algorithm}
  \KwIn{An even integer $2n \geq 6$, representing the required number of obstacle vertices.\\
  The width $w$ and height $h$ of the bounding box of the obstacle.\\
  The obstacle type (L-obstacle, T-obstacle, staircase obstacle, partial staircase or general obstacle).}
  \KwOut{A rectilinearly-convex obstacle of the desired type with $2n$ vertices.}
  // Let $n_b$ and $n_t$ denote the number of horizontal boundary edges visible from below and above the obstacle, respectively, and let $i_b$ and $i_t$ denote the labels of the bottom and top extreme edges respectively (where the labels are increasing integers, from left to right, starting at 1).\\
  \nl \If{$\omega$ is an L-obstacle} {
        \nl $n_b := 1$; $n_t := n - 1$; $i_b := 1$; $i_t := 1$.\\
  } \nl \ElseIf{$\omega$ is a T-obstacle} {
        \nl $n_b := n - 1$; $n_t := 1$; $i_b$ is selected at random from $\{2,\ldots,n_b-1\}$; $i_t := 1$.
    } \nl \ElseIf{$\omega$ is a staircase obstacle} {
        \nl $n_b$ is selected at random from $\{2,\ldots,n-2\}$; $n_t := n - n_b$; $i_b := 1$; $i_t := n_t$.
    } \nl \ElseIf{$\omega$ is a partial staircase} {
        \nl $n_b$ is selected at random from $\{3,\ldots,n-2\}$; $n_t := n - n_b$; $i_b$ is selected at random from $\{2,\ldots,n_b-1\}$; $i_t := n_t$.
    } \nl \Else{
            $n_b$ is selected at random from $\{3,\ldots,n-3\}$; $n_t := n - n_b$; $i_b$ is selected at random from $\{2,\ldots,n_b-1\}$; $i_t$ is selected at random from $\{2,\ldots,n_t-1\}$.
    }
    \nl // Determine the number of vertical boundary edges $n_l$ and $n_r$ visible from the left and right of the obstacle respectively.\\
    \nl $n_l := (i_b - 1) + (i_t - 1) + 1$; $n_r := n - n_l$.\\
    \nl Randomly partition $h$ into $n_l$ segments corresponding to the vertical edges visible to the left of the obstacle and $n_r$ segments corresponding to the vertical edges visible to the right of the obstacle.\\
    \nl Randomly partition $w$ into $n_b$ segments corresponding to the horizontal edges visible to the bottom of the obstacle and $n_t$ segments corresponding to the horizontal edges visible to the top of the obstacle.\\
    \nl Starting at the left vertex of the lower extreme edge, move around the perimeter counterclockwise and alternately add horizontal and vertical edges to the boundary of the obstacle using the above partitions.\\
    \nl If the boundary of intersects itself during construction, re-sample the current partition accordingly (see for instance the example in Figure~\ref{fig-construct-random-rcpoly} (c)).
  \caption{\textsc{RandomRectilinearlyConvexObstacle}}
  \label{RandomRectilinearlyConvexObstacle}
\end{algorithm}

\begin{figure}[h]
\centering
\includegraphics[width=\textwidth]{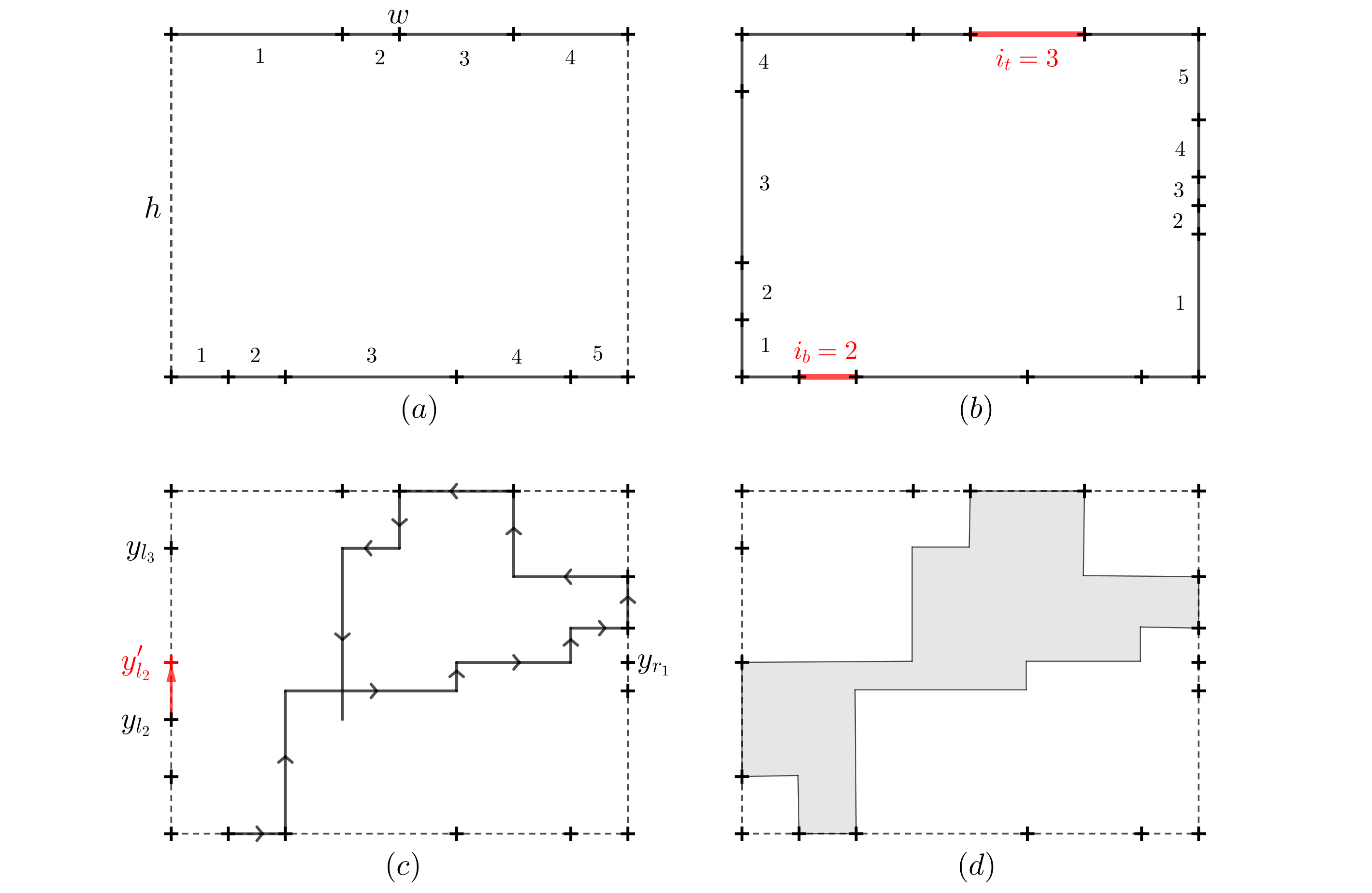}
\caption{Computing a random rectilinearly-convex polygon using Algorithm~\ref{RandomRectilinearlyConvexObstacle}.}\label{fig-construct-random-rcpoly}
\end{figure}

Figure~\ref{fig-construct-random-rcpoly} shows the construction of a random rectilinearly-convex polygon with $2n = 18$ vertices, width $w$ and height $h$. In (a), $n_b = 5$ is randomly computed, and therefore $n_t = 9 - 5 = 4$. The width is randomly partitioned into segments of varying length for the bottom and top sides of the bounding box. In (b), the bottom and top extreme edges are randomly selected to be $i_b = 2$ and $i_t = 3$ respectively. The number of segments on the left and right sides of the bounding box are determined as follows: $n_l = (i_b - 1) + (i_t - 1) + 1 = 4$ and $n_r = (n_b - i_b) + (n_t - i_t) + 1 = 5$. The left and right sides of the bounding box are each partitioned into an appropriate number of segments of varying length. In (c), a self-intersection is detected, and the point $y_{l_2}$ is moved to $y_{l_2}'$ where $y_{l_2}'$ is randomly sampled at a location between $y_{r_1}$ and $y_{l_3}$. The final obstacle is shown in (d). Note that if the obstacle type is not specified, most of the obstacles generated by Algorithm~\ref{RandomRectilinearlyConvexObstacle} will be general obstacles.

\subsection{Computational results}

Table~\ref{table-experiments} shows the number of edges in a minimum skeleton for random rectilinearly-convex obstacles with 100 to 1000 vertices. For comparison, the results for strictly-connected skeletons are shown in parentheses. The values stated are for 100 seeds. Experiments were not run for the trivial cases when the obstacle is a rectangle, an L-obstacle or a T-obstacle (in which cases the number of edges in a minimum skeleton is either one or two). The results indicate that in the worst case (for staircase obstacles) a minimum skeleton for a random rectilinearly-convex obstacle with 1000 vertices has a small number of edges compared to the number of edges in the boundary of the original obstacle. Staircase obstacles generally require more skeleton edges compared to partial staircase obstacles and general obstacles, since staircase obstacles have a higher likelihood of having a `long skinny' structure due to their pairs of extreme edges lying at opposite ends of the obstacle. General obstacles generally require fewer skeleton edges due to their having a higher likelihood of having a cross, opposite extreme visibility edge or adjacent extreme visibility edge. For a given partial staircase or general obstacle, as the number of vertices increases, the obstacle grows into the convex hull of its extreme edges. As a consequence, partial staircases and general obstacles with many vertices generally have only a small number of skeleton edges. The results indicate that, in general, strictly-connected skeletons do not have significantly more edges compared to weakly-connected skeletons.

\begin{table}[ht]
  \centering
  {\tabcolsep=0pt\def\arraystretch{1}
\noindent\begin{tabularx}{\columnwidth}{@{\extracolsep{\stretch{2}}}*{10}{r}@{}}\toprule
\centering
   & \multicolumn{3}{c}{Staircase} & \multicolumn{3}{c}{Partial staircase} & \multicolumn{3}{c}{General obstacle}\tabularnewline
   \cmidrule(l){2-4} \cmidrule(l){5-7} \cmidrule(l){8-10}
  Vertices & Min & Median & Max & Min & Median & Max & Min & Median & Max \tabularnewline
  100  & 2 (2) & 7 (8) & 20 (22) & 2 (2) & 3 (4) & 18 (19) & 2 (2) & 3 (3) & 11 (11) \tabularnewline
  200  & 1 (1) & 13 (14) & 37 (39) & 2 (2) & 3 (4) & 17 (18) & 2 (2) & 3 (3) & 9 (9) \tabularnewline
  300  & 2 (2) & 13 (14) & 48 (49) & 2 (2) & 3 (4) & 40 (42) & 2 (2) & 3 (3) & 15 (16) \tabularnewline
  400  & 2 (2) & 19 (21) & 63 (64) & 2 (2) & 3 (4) & 56 (57) & 2 (2) & 3 (3) & 9 (10) \tabularnewline
  500  & 2 (2) & 20 (23) & 79 (80) & 2 (2) & 3 (4) & 45 (46) & 2 (2) & 3 (3) & 38 (42) \tabularnewline
  600  & 2 (2) & 22 (26) & 80 (86) & 2 (2) & 3 (4) & 66 (67) & 2 (2) & 3 (3) & 17 (18) \tabularnewline
  700  & 2 (2) & 25 (29) & 101 (103) & 2 (2) & 3 (4) & 30 (31) & 2 (2) & 3 (3) & 6 (7) \tabularnewline
  800  & 2 (2) & 24 (27) & 86 (94) & 2 (2) & 3 (4) & 35 (36) & 2 (2) & 3 (3) & 5 (5) \tabularnewline
  900  & 2 (2) & 24 (27) & 101 (111) & 2 (2) & 3 (4) & 77 (84) & 2 (2) & 3 (3) & 7 (8) \tabularnewline
  1000 & 2 (2) & 33 (36) & 142 (152) & 2 (2) & 4 (5) & 42 (45) & 2 (2) & 3 (3) & 8 (8) \tabularnewline\bottomrule
\end{tabularx}}
  \caption{Number of edges in a minimum skeleton by obstacle type and varying the number of obstacle vertices. Results for strictly-connected skeletons are shown in parentheses. Values stated are for 100 seeds.
  }\label{table-experiments}
\end{table}

Preliminary tests were undertaken for randomly-generated obstacles constructed using a modification of Algorithm~\ref{RandomRectilinearlyConvexObstacle} in which the obstacle edge lengths were randomly determined according to an exponential distribution rather than a uniform distribution (the resulting obstacles have small numbers of edges that are relatively large compared to the other edges), however these tests did not yield significantly different results.

\begin{figure}[h]
\centering
\includegraphics[width=\textwidth]{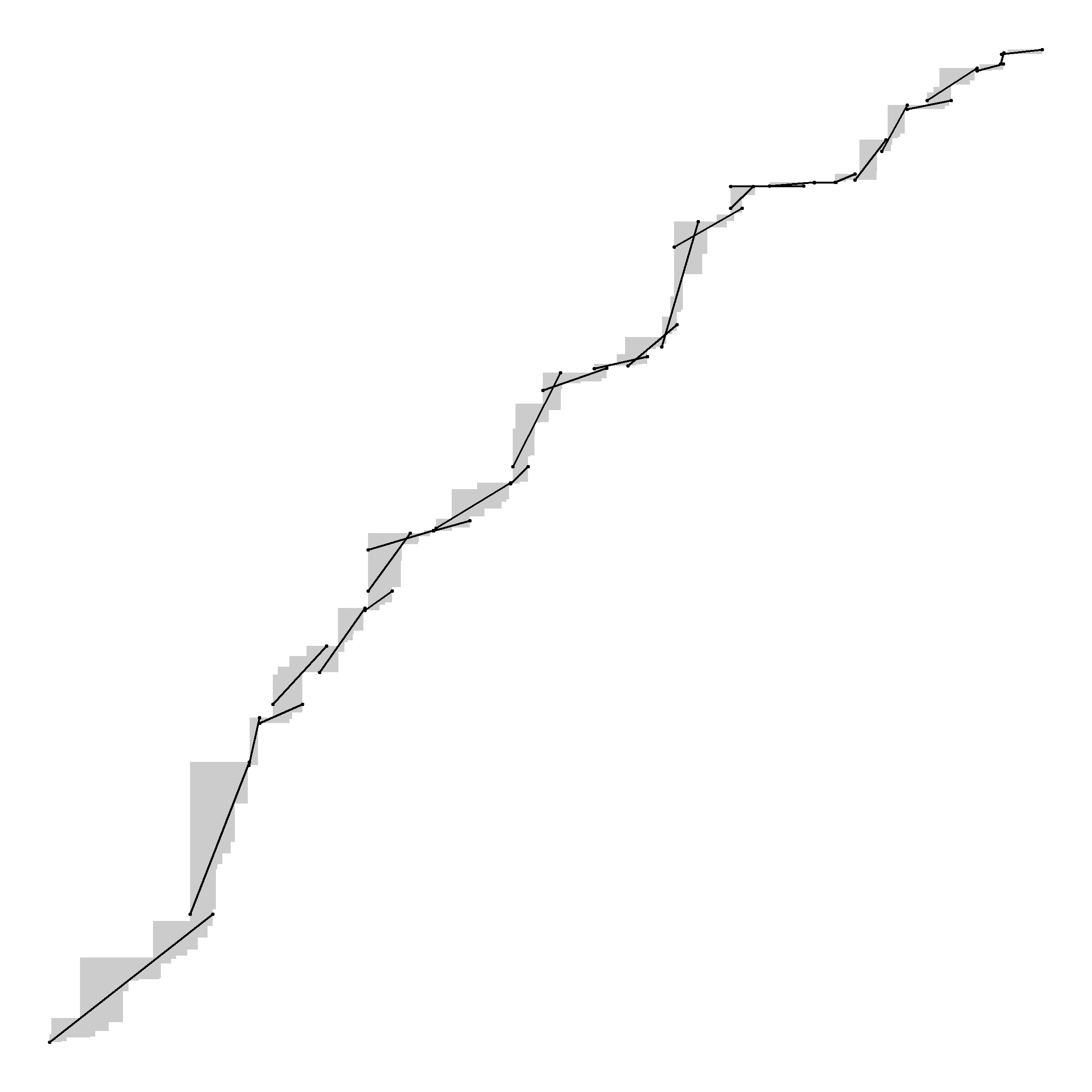}
\caption{A minimum skeleton for a staircase obstacle with 300 vertices. The skeleton has 28 edges.}\label{staircase_n300_s19}
\end{figure}

An example of a minimum skeleton for a staircase obstacle with 300 vertices is shown in Figure~\ref{staircase_n300_s19}.

\section{Conclusion}

We have introduced the concept of an obstacle \emph{skeleton} which is a set of line segments inside an obstacle $\omega$ that can be used in place of $\omega$ when performing intersection tests for obstacle-avoiding shortest network problems in the plane. A skeleton can have significantly fewer line segments compared to the number of line segments in the boundary of the original obstacle, and therefore performing intersection tests on a skeleton (rather than the original obstacle) can potentially significantly reduce the CPU time required by algorithms for computing shortest obstacle-avoiding networks. We have provided an exact $O(n^2)$ algorithm for computing minimum skeletons for obstacles in the rectilinear plane that are rectilinearly-convex (obstacles whose edges are either horizontal or vertical and for which any two points in the obstacle have a shortest rectilinear path that is entirely inside the obstacle), in the context of the obstacle-avoiding rectilinear Steiner tree problem. We have shown that skeletons for rectilinearly-convex obstacles have at most half as many edges compared to the boundaries of the original obstacles, and that, in practice, the number of edges in a minimum skeleton is significantly smaller. Future work will look at heuristics and exact algorithms for computing skeletons for convex and non-convex obstacles in other fixed orientation metrics as well as the Euclidean metric.

\section{Acknowledgements}

We thank Martin Zachariasen for interesting discussions and feedback over the course of developing this work. This work was supported by an Australian Research Council Discovery Grant.

\bibliographystyle{spmpsci_unsrt}
\bibliography{references}

\end{document}